\def\version{\today}

\documentclass[12pt]{article}

\usepackage{amssymb}
\usepackage{dsfont}
\usepackage{graphicx}
\usepackage{bbm}
\usepackage[usenames,dvipsnames]{color}
\usepackage{amsthm}
\usepackage{amsmath}
\usepackage{index}
\usepackage{enumerate}
\usepackage{ifpdf}
\usepackage[pdftex]{thumbpdf}       
\graphicspath{{pictures/}}

\usepackage{eepic}
\usepackage{multirow}

\ifpdf
\usepackage[a4paper,pdftex]{hyperref}

  \hypersetup{%
  pdfauthor={Remco van der Hofstad},%
  pdftitle={Random Graphs and Complex Networks},%
  pdfsubject={Monograph},%
  pdfkeywords={},%
  pdfstartview={FitH},%
  pdfpagelayout={OneColumn},%
  hyperindex=true,%
  bookmarks=true,%
  bookmarksnumbered=true,%
  bookmarksopen=false,%
  unicode=true,%
  breaklinks=true,%
  colorlinks=true,%
  citecolor=MidnightBlue,%
  linkcolor=BlueViolet,%
  plainpages=false,
}%
\else
  \newcommand\phantomsection\relax
  \newcommand{\url}[1]{#1}  
  \newcommand{\href}[2]{#2} 
  \usepackage{nohyperref}
\fi

\newcommand{\e}{{\mathrm e}}
\renewcommand{\i}{{\mathrm i}}

\theoremstyle{plain}
\newtheorem{theorem}{Theorem}[section]
\newtheorem{corollary}[theorem]{Corollary}
\newtheorem{lemma}[theorem]{Lemma}

\newtheorem{proposition}[theorem]{Proposition}

\newtheorem{op}[theorem]{Open problem}

\theoremstyle{remark}

\newtheoremstyle{maria}{}{}{}{0cm}{\itshape}{}{\newline}{}
\theoremstyle{maria}

\numberwithin{equation}{section}
\numberwithin{theorem}{section}

\usepackage{fancyhdr}

\newcommand{\bfwit}{\boldsymbol{w}}
\newcommand{\bfdit}{\boldsymbol{d}}
\newcommand{\bfpit}{\boldsymbol{p}}

\newcommand{\Nbold} {{\mathbb N}}

\newcommand{\Rbold} {{\mathbb R}}

\newcommand{\Ccal}   {\mathcal{C}}

\renewcommand{\Re}{{\rm Re}}
\renewcommand{\Im}{{\rm Im}}

\newcommand{\ERnp}{{\rm ER}_n(p)}

\newcommand{\GRGnw}{{\rm GRG}_n(\bfwit)}

\newcommand{\CLnw}{{\rm CL}_n(\bfwit)}

\newcommand{\NRnw}{{\rm NR}_n(\bfwit)}

\newcommand{\Ver}{V_n}

\newcommand{\BIN}{{\rm BIN}}

\newcommand{\Poi}{{\rm Poi}}

\newcommand {\convp}{\stackrel{\sss {\mathbb P}}{\longrightarrow}}

\newcommand {\SD}{\preceq}

\renewcommand{\op}{o_{\prob}}

\def\1{{\mathchoice {1\mskip-4mu\mathrm l}      
{1\mskip-4mu\mathrm l}
{1\mskip-4.5mu\mathrm l} {1\mskip-5mu\mathrm l}}}
\newcommand{\indic}[1]{\1_{\{#1\}}}
\newcommand{\indicwo}[1]{\1_{#1}}

\newcommand{\expec}{{\mathbb{E}}}
\newcommand{\prob}{{\mathbb{P}}}

\newcommand{\Var}{{\rm Var}}

\newcommand{\R}{\Rbold}

\newcommand{\N}{\Nbold}

\renewcommand{\e}{{\mathrm e}}
\newcommand{\tildeX}{\widetilde{X}}
\newcommand{\tildeS}{\widetilde{S}}
\newcommand{\tildeT}{\widetilde{T}}

\newcommand{\sss}   { \scriptscriptstyle }

\newcommand {\sX}{{\sss X}}
\newcommand {\sY}{{\sss Y}}

\newcommand {\vep}{\varepsilon}

\newcommand{\conn}{\longleftrightarrow}

\newcommand{\nc}        { \conn {\hspace{-2.8ex} /} \hspace{1.8ex}   }

\newcommand{\Cmax}{{\mathcal C}_{\rm max}}
\newcommand{\cluster}{{\mathcal C}}

\def\eqalign#1\enalign{
    \begin{align}#1\end{align}
    }
\newcommand{\lbeq}[1]  {\label{e:#1}}
\newcommand{\refeq}[1] {\eqref{e:#1}}
\newcommand{\eq}{\begin{equation}}
\newcommand{\en}{\end{equation}}
\newcommand{\ben}{\begin{enumerate}}
\newcommand{\een}{\end{enumerate}}
\newcommand{\eqn}[1]{\begin{equation} #1 \end{equation}}
\newcommand{\eqan}[1]{\eqalign #1 \enalign}

\newcommand{\nn}{\nonumber}

\renewcommand{\to}{\rightarrow}

\title  {
        Critical behavior in inhomogeneous random graphs
        }

\author{
{Remco van der Hofstad}
\thanks{Department of Mathematics and Computer Science,
Eindhoven University of Technology, P.O.\ Box  513,
5600 MB Eindhoven, The Netherlands.
{\tt rhofstad@win.tue.nl}}
}

\date\version

\begin{document}

\maketitle

\begin{abstract}
We study the critical behavior of inhomogeneous
random graphs where edges are present independently but with unequal
edge occupation probabilities. The edge probabilities are moderated by
\emph{vertex weights}, and are such that the degree of vertex $i$
is close in distribution to a Poisson random variable with parameter
$w_i$, where $w_i$ denotes the weight of vertex $i$. We choose the weights
such that the weight of a uniformly chosen vertex converges in distribution
to a limiting random variable $W$, in which case the proportion of
vertices with degree $k$ is close to the probability that a Poisson
random variable with \emph{random} parameter $W$ takes the value $k$.
We pay special attention to the \emph{power-law case},
in which $\prob(W\geq k)$ is proportional to $k^{-(\tau-1)}$
for some power-law exponent $\tau>3$, a property which is
then inherited by the asymptotic degree distribution.

We show that the critical behavior depends sensitively on the
properties of the asymptotic degree distribution moderated by the
asymptotic weight distribution $W$. Indeed, when $\prob(W\geq k)
\leq ck^{-(\tau-1)}$ for all $k\geq 1$ and some $\tau>4$ and $c>0$,
the largest critical connected component in a graph of size $n$ is of
order $n^{2/3}$, as on the Erd\H{o}s-R\'enyi random graph.
When, instead, $\prob(W\geq k)=ck^{-(\tau-1)}(1+o(1))$ for $k$ large
and some $\tau\in (3,4)$ and $c>0$, the largest
critical connected component is of the much smaller
order $n^{(\tau-2)/(\tau-1)}.$

\end{abstract}

\section{Introduction and results}
\label{sec-intro}

We study the critical behavior of inhomogeneous
random graphs, where edges are present independently but with unequal
edge occupation probabilities. Such inhomogeneous random graphs
were studied in substantial detail in the seminal paper
by Bollob\'as, Janson and Riordan \cite{BolJanRio07},
where various results have been proved,
including their critical value by studying the
connected component sizes in the super- and subcritical regimes.

In this paper, we study the \emph{critical behavior}
of such random graphs, and show that this critical
behavior depends sensitively on the asymptotic properties
of their degree sequence, i.e., the asymptotic proportion
of vertices with degree $k$ for each $k\geq 1$.
Our results show that the critical behavior of our inhomogeneous random graphs
admits a transition when the third moment of the degrees turns from finite to infinite.

\subsection{Inhomogeneous random graphs: the rank-1 case}
\label{sec-IRG}
In this section, we introduce the random graph model that we shall investigate.
In our models, $\bfwit=(w_j)_{j\in [n]}$ are vertex weights,
and $\ell_n$ is the total weight of all vertices given by $\ell_n=\sum_{j=1}^n w_j$.
We shall mainly work with the Poisson random graph or {\it Norros-Reittu random graph}
\cite{NorRei06}, which we denote by $\NRnw$. In the $\NRnw$,
the edge probabilities are given by
    \eqn{
    \lbeq{edgeprobNR}
    p_{ij}^{\sss({\rm NR})}=1-\mathrm{e}^{-w_iw_j/\ell_n}.
    }
More precisely, $p_{ij}^{\sss({\rm NR})}$ is the probability that edge
$ij$ is present or \emph{occupied}, for $1\leq i<j\leq n$, and different edges are
independent. In Section \ref{sec-asy-eq},
we shall extend our results to graphs where the edge probabilities are
either $p_{ij}=\max\{w_iw_j/\ell_n, 1\}$ (as studied by Chung and Lu
in \cite{ChuLu02b,ChuLu03, ChuLu06c, ChuLu06}) or $p_{ij}=w_iw_j/(\ell_n+w_iw_j)$
(as studied by Britton, Deijfen and Martin-L\"of in \cite{BriDeiMar-Lof05}).
See \cite[Section 16.4]{BolJanRio07} for a detailed discussion of the relation
between the general inhomogeneous random graph and the models studied here,
which are called \emph{rank-1 inhomogeneous random graphs} in \cite{BolJanRio07}.

Naturally, the graph structure depends sensitively on
the empirical properties of the weights, which we shall now
introduce. Let $F$ be a distribution
function, and define
    \eqn{
    \lbeq{choicewi}
    w_{j} = [1-F]^{-1}(j/n),
    }
where $[1-F]^{-1}$ is the \emph{generalized inverse} of
$1-F$ defined, for $u\in (0,1)$, by
    \eqn{
    \lbeq{invverd}
    [1-F]^{-1}(u)=\inf \{ s: [1-F](s)\leq u\}.
    }
By convention, we set $[1-F]^{-1}(1)=0$.

In the setting in \refeq{edgeprobNR} and \refeq{choicewi}, by
\cite[Theorem 3.13]{BolJanRio07}, the
number of vertices with degree $k$, which we denote by $N_k$,
satisfies, with $W$ having distribution function $F$
appearing in \refeq{choicewi},
    \eqn{
    \lbeq{N-k-conv}
    N_k/n\convp f_k\equiv \expec\Big[\mathrm{e}^{-W} \frac{W^k}{k!}\Big], \qquad k\geq 0,
    }
where $\convp$ denotes convergence in probability.
We recognize the limiting distribution $\{f_k\}_{k=1}^{\infty}$ as a
so-called \emph{mixed Poisson distribution with mixing distribution $F$},
i.e., conditionally on $W=w$, the distribution is
Poisson with mean $w$. Since a Poisson random variable with a large parameter
is highly concentrated around that parameter, it is intuitively clear that the number of vertices with
degree larger than $k$ is, for large $k$, quite close to $n[1-F(k)]$.
In particular, for $a>0$, $\sum_k k^a f_k<\infty$ precisely when
$\expec[W^a]<\infty$.

In our setting, there exists a giant component containing a
positive proportion of the vertices precisely when
$\nu>1$, where we define
    \eqn{
    \lbeq{nu-def}
    \nu=\frac{\expec[W^2]}{\expec[W]}.
    }
As we explain in more detail in Section \ref{sec-disc}, we shall see that $\nu$ arises
as the mean of the \emph{size-biased distribution} of $W$, which, in turn,
arises as the mean offspring in a branching process approximation of the
exploration of the connected component of a vertex. More precisely, if $\nu>1$,
then the largest connected component has $n\zeta(1+\op(1))$ vertices, while if
$\nu\leq 1$, the largest connected component has $\op(n)$ vertices.
Here we write that $X_n=\op(b_n)$ for some sequence $b_n$,
when $X_n/b_n$ converges to zero in probability.
See, e.g., \cite[Theorem 3.1 and Section 16.4]{BolJanRio07}
and \cite{ChuLu02b,ChuLu06,NorRei06}. When $\nu>1$, the rank-1
inhomogeneous random graph is called {\it supercritical}, when $\nu=1$
it is called {\it critical}, and when $\nu<1$, it is called {\it subcritical}.
The aim of this paper is to study the size of the largest connected components
in the critical case.

\subsection{Results}
\label{sec-results}
Before we can state our results, we introduce some notation.
We write $[n]=\{1, \ldots, n\}$ for the set of vertices.
For two vertices $s, t\in [n]$, we write $s\conn t$ when
there exists a path of occupied edges
connecting $s$ and $t$. By convention, $v\conn v$.
For $v\in [n]$, we denote the {\it cluster of $v$} by
$\cluster(v)=\big\{x\in [n]\colon v\conn x\big\}$.
We denote the size of $\cluster(v)$ by $|\cluster(v)|$, and define
the \emph{largest connected component} by
$|\Cmax|=\max\{|\cluster(v)|\colon v\in [n]\}$.
Our main results are:

\begin{theorem}[Largest critical cluster for $\tau>4$]
\label{main-clusterdistr-tau>4}
Fix $\NRnw$ with $\bfwit=(w_j)_{j\in [n]}$ as in \refeq{choicewi}, and assume that
the distribution function $F$ in \refeq{choicewi} satisfies $\nu=1$.
Suppose there exists a $\tau>4$ and a constant $c_{\sss F}>0$ such that,
for all large enough $x\geq 0$,
    \eqn{
    \lbeq{F-bound-tau>4}
    1-F(x)\leq c_{\sss F}x^{-(\tau-1)}.
    }
Let $\tilde{\bfwit}=(\tilde{w}_j)_{j\in [n]}$ be defined by
    \eqn{
    \lbeq{tilde-wi-def}
    \tilde{w}_j=(1+\vep_n)w_j,
    }
and fix $\vep_n$ such that $|\vep_n|\leq \Lambda n^{-1/3}$ for some $\Lambda>0$.
Then there exists a constant $b=b(\Lambda)>0$ such that for all $\omega>1$
and for $n$ sufficiently large, ${\rm NR}_n(\tilde{\bfwit})$ satisfies
    \eqn{
    \lbeq{crit-window-tau>4}
    \prob\Big(\omega^{-1} n^{2/3}\leq |\Cmax|\leq \omega n^{2/3}
    \Big)
    \geq 1-\frac{b}{\omega}.
    }
\end{theorem}

\begin{theorem}[Largest critical cluster for $\tau\in (3,4)$]
\label{main-clusterdistr-tau(3,4)}
Fix $\NRnw$ with $\bfwit=(w_j)_{j\in [n]}$ as in \refeq{choicewi}, and assume that
the distribution function $F$ in \refeq{choicewi} satisfies $\nu=1$.
Suppose that there exists a $\tau\in (3,4)$ and a constant
$0<c_{\sss F}<\infty$ such that
    \eqn{
    \lbeq{F-bound-tau(3,4)b}
    \lim_{x\rightarrow \infty} x^{\tau-1}[1-F(x)]= c_{\sss F}.
    }
Fix $\vep_n$ such that $|\vep_n|\leq \Lambda n^{-(\tau-3)/(\tau-1)}$ for some $\Lambda>0$.
Then there exists a constant $b=b(\Lambda)>0$ such that for all $\omega>1$
and for $n$ sufficiently large, ${\rm NR}_n(\tilde{\bfwit})$, with $\tilde{\bfwit}$
defined as in \refeq{tilde-wi-def}, satisfies
    \eqn{
    \lbeq{crit-window-tau(3,4)b}
    \prob\Big(\omega^{-1} n^{(\tau-2)/(\tau-1)}\leq |\Cmax|
    \leq \omega n^{(\tau-2)/(\tau-1)}\Big)\geq 1-\frac{b}{\omega}.
    }
\end{theorem}

\subsection{Discussion and related results}
\label{sec-disc}
In this section, we discuss our results and the relevant results
in the literature. We start by introducing some notation used throughout this paper.
We write $X\sim \Poi(\lambda)$ to denote that $X$ has a Poisson distribution with
(possibly random) parameter $\lambda$, and $a_n=\Theta(b_n)$ if there exist positive constants
$c$ and $C$, such that, for all $n$, we have $cb_n\le a_n\le Cb_n$.

\paragraph{Branching process approximation.}
The main tool used in this paper is the comparison of clusters
to \emph{branching processes.} Let $\Ver$ be a vertex chosen uniformly
from $[n]$. Then, the number of neighbors of $\Ver$ is close to
$\Poi(W_n)$, where $W_n=w_{\Ver}$ is the (random) weight of
$\Ver$. In the setting of \refeq{choicewi}, we shall see that
$W_n$ converges in distribution to a random variable $W$
having distribution function $F$, which explains that a uniformly chosen vertex has
a degree that is close to $\Poi(W)$ (recall \refeq{N-k-conv}).
As described in more detail in Section \ref{sec-MPBP}, we can
describe the set of vertices to which $\Ver$ is connected
by associating a \emph{random mark} to each of the $\Poi(W_n)$
values, where the mark equals $i\in [n]$ with probability
$w_i/\ell_n$. Then, the set of neighbors of $\Ver$ equals
the set of marks chosen. Further, the distribution of the
\emph{degree} of a neighbor of $\Ver$ is close to $X\sim \Poi(w_{M})$, where
$M$ is the mark associated to the neighbor, and
the degrees of different neighbors of $\Ver$ are close to an i.i.d.\ sequence.
Thus, the cluster exploration is close to a
\emph{branching process with offspring distribution}
$X\sim \Poi(w_{\sss M})$. It is not hard to see that
$\Poi(w_{\sss M})$ converges in distribution to
$\Poi(W^*)$, where, for a non-negative random variable $X$ with $\expec[X]>0$,
we let $X^*$ denote its \emph{size-biased distribution} given by
    \eqn{
    \lbeq{size-biased-distr}
    \prob(X^*\leq x)=\frac{\expec[X\indic{X\leq x}]}{\expec[X]}.
    }
Thus, $\Poi(W^*)$ has finite variance when $W$ has a finite third moment.
For details, see Proposition \ref{prop-N NR kop},
where this connection is made explicit.
We denote the \emph{mean offspring} of the branching process by
    \eqn{
    \lbeq{nu-n-def}
    \nu_n=\expec[\Poi(w_{M})]=\frac{\sum_{j=1}^n w_j^2}{\sum_{j=1}^n w_j}.
    }
In the setting of \refeq{choicewi}, we shall see that $\nu_n\rightarrow \nu$,
where $\nu$ is defined by \refeq{nu-def} (see Corollary \ref{cor-varphi-n}(b) below).
Therefore, the resulting branching process is \emph{critical}
precisely when $\nu=1$. Observe that the offspring
$\Poi(W^*)$ of this branching process has
\emph{finite variance} when $\tau>4$, but not when $\tau\in (3,4)$.
We now make use of the relation to branching processes to
connect the subcritical and supercritical regimes to the critical one.

\paragraph{Connecting the subcritical and supercritical regimes
to the critical one.}
We first give a heuristic explanation for the critical behavior
of $n^{2/3}$ appearing in Theorem \ref{main-clusterdistr-tau>4}.
Let $\vep_n=\nu_n-1$ and $\tau>4$. By the branching process approximation,
the largest connected component has size $\rho_n n(1+o(1))$ when $\vep_n>0$,
where $\rho_n$ is the survival probability of
the branching process approximation to the cluster.
Now, $\rho_n$ is of the order $\vep_n$ when $\tau>4$,
since the corresponding branching process has finite variance in this case.
On the other hand, the largest subcritical cluster is
$\Theta(\vep_n^{-2}\log{(n\vep_n^3)})$ when $\vep_n<0$,
since, for branching processes with finite mean, the probability that
the total progeny exceeds $k$ is approximately equal to
$\Theta(1/\sqrt{k})\e^{-\Theta(k \vep_n^2)}$. This suggests that
the critical behavior arises precisely when
$\vep^{-2}_n=n\vep_n$, i.e., when $\vep_n=n^{-1/3}$,
and in this case, the largest connected component is
$\vep_n n=n^{2/3}$ as in Theorem
\ref{main-clusterdistr-tau>4}.

We next extend this heuristic to the case $\tau\in (3,4)$, for which
the picture changes completely. The results by Janson in \cite{Jans08b}
suggest that the largest subcritical cluster
is like $w_1/(1-\nu)=\Theta(n^{1/(\tau-1)}/|\vep_n|)$ when $\nu_n=1+\vep_n$
and $\vep_n<0$. We note that
\cite{Jans08b} only proves this when $\nu<1$ is \emph{fixed},
but we conjecture that it extends to all subcritical $\nu$.
In the supercritical regime, instead, the largest connected component
should be like $n \rho_n$, where $\rho_n$ is the
survival probability of the (infinite variance) branching process
approximation of the cluster. A straightforward computation shows that,
when $\vep_n>0$ and $\vep_n=o(1)$, we have
$\rho_n \sim \vep_n^{1/(\tau-3)}$ (see Lemma \ref{lem-surv-prob}
below). Thus, this suggests that the critical behavior should now
be characterized instead by taking $n^{1/(\tau-1)}/\vep_n=\vep_n^{1/(\tau-3)} n$,
which is $\vep_n=n^{-(\tau-3)/(\tau-1)}$.
In this case, the largest critical cluster should be of the order
$\vep_n^{1/(\tau-3)} n\sim n^{(\tau-2)/(\tau-1)}$,
as in Theorem \ref{main-clusterdistr-tau(3,4)}.
This suggests that in both cases, the subcritical and supercritical
regimes connect up nicely. In order to make these statements precise,
and thus showing that Theorems
\ref{main-clusterdistr-tau>4}--\ref{main-clusterdistr-tau(3,4)} really
deal with all the `critical weights', we would need to show that
when $|\vep_n|$ is much larger than $n^{-1/3}$ and $n^{-(\tau-3)/(\tau-1)}$,
respectively, the above heuristic bounds on $|\Cmax|$, for both
the super- and subcritical regimes, are precise.


\paragraph{The scaling limit of cluster sizes for $\tau>4$.}
A special case of Theorem \ref{main-clusterdistr-tau>4} is
the critical behavior for the Erd\H{o}s-R\'enyi random graph (ERRG), where bounds as in
\refeq{crit-window-tau>4} have a long history (see e.g., \cite{ErdRen60}, as well as
\cite{Boll84b,JanKnuLucPit93,LucPitWie94,Pitt01} and the monographs \cite{Boll01,JanLucRuc00}
for the most detailed results). The ERRG corresponds to
taking $w_j=c$ for all $j\in [n]$, and then $\nu$ in \refeq{nu-def} equals $c$.
Therefore, criticality corresponds to $w_j=1$ for all $j\in [n]$.
For the ERRG there is a tremendous amount of work on
the question for which values of $p$, similar critical behavior
is observed as for the critical value $p=1/n$
\cite{Aldo97, Boll84b, ErdRen60, JanKnuLucPit93, Lucz90a, LucPitWie94}.
Indeed, when we take $p=(1+\lambda n^{-1/3})/n$, the largest cluster
has size $\Theta(n^{2/3})$ for every fixed $\lambda\in \R$, but
it is $\op(n^{2/3})$ when $\lambda\rightarrow -\infty$, and has size $\gg n^{2/3}$
when $\lambda\gg 1$. Therefore, the values $p$ satisfying $p=(1+\lambda n^{-1/3})/n$
for some $\lambda\in {\mathbb R}$ are sometimes called the \emph{critical window.}

Aldous \cite{Aldo97} proves that the vector of ordered
cluster sizes of the ERRG weakly converges to
a limiting process, which can be characterized as the excursions
of a standard Brownian motion with a parabolic drift,
ordered in their sizes. A less well-known extension by Aldous
\cite{Aldo97} can be found in \cite[Prop. 4]{Aldo97},
where an inhomogeneous random graph is studied in which
there is an edge between $i$ and $j$ with probability $1-\e^{-x_ix_j}$ for some
vertex weights $(x_i)_{i\in [n]}$ and different edges are
independent. This corresponds to our setting when we take
$x_i=w_i/\sqrt{\ell_n}$. Aldous shows in \cite[Prop. 4]{Aldo97} that
the \emph{ordered cluster weights} weakly converge to a limit closely
related to that of the ERRG, where the weight of a set of vertices $\Ccal$
equals $\sum_{c\in \Ccal} x_c$.  Since the completion of the first version of this paper,
in fact the weak convergence of the ordered cluster sizes has
been proved independently and almost at the same time
in \cite{Turo09, BhaHofLee09a}, using related means as in
\cite{Aldo97}, and under the slightly weaker
condition that $\expec[W^3]<\infty$. We have included the proof of
Theorem \ref{main-clusterdistr-tau>4} as this proof
follows the same lines as the proof of the novel result in
Theorem \ref{main-clusterdistr-tau(3,4)}, and the proof nicely
elucidates the place where the restriction $\tau>4$ is used.
It is not hard to see that our proofs in fact carry
over to the situation where $\expec[W^3]<\infty$, but we
refrain from doing so for simplicity.

\paragraph{The scaling limit of cluster sizes for $\tau\in (3,4)$.}
When $\tau\in (3,4)$, large parts of the above discussion remain valid,
however, the variance of $\Poi(W^*)$ arising in the exploration process is
{\it infinite}. Therefore, the critical nature of the total progeny
of the branching process approximation is rather different,
which is reflected in different critical behavior.
Since the completion of the first version of this paper,
in fact the weak convergence of the ordered cluster sizes has
been proved in \cite{BhaHofLee09b}.
The proof relies on the fact that the cluster exploration
can be described by a \emph{thinned L\'evy process} having rather
interesting behavior. In the proof in \cite{BhaHofLee09b},
the results derived in this paper, in particular Propositions
\ref{prop-main} and \ref{thm-UBclust-size} below, play a crucial role.

In our results, we have assumed the precise power-law form
of $1-F$ in \refeq{F-bound-tau(3,4)b}. A heuristic computation shows
that when $u\mapsto [1-F]^{-1}(u)=u^{-1/(\tau-1)}\ell(u)$ for some
function $u\mapsto \ell(u)$ slowly
varying at $u=0$, then the size of $|\Cmax|$ becomes $n^{(\tau-2)/(\tau-1)}/\ell(1/n)$,
and the width of the critical window becomes $n^{(\tau-3)/(\tau-1)} \ell(1/n)^2$.
This also sheds light on the \emph{critical cases} $\tau=3$ and $\tau=4$. Indeed,
when $\tau=3$ and $u\mapsto \ell(u)$ is such that
$\expec[W^2]<\infty$, we predict that the above applies. If $\tau=4$ and
$u\mapsto \ell(u)$ is such that $\expec[W^3]<\infty$, then we predict that
$|\Cmax|$ is of order $n^{2/3}$ as in Theorem \ref{main-clusterdistr-tau>4},
while if $\tau=4$ and $u\mapsto \ell(u)$ is such that $\expec[W^3]=\infty$, then we predict that
$|\Cmax|$ is of order $n^{2/3}/\ell(1/n)$, instead. The predictions for the critical window
are accordingly. In our proofs, the presence of a slowly varying function should enter
in Propositions \ref{prop-main}--\ref{thm-UBclust-size} and Lemma \ref{lem-surv-prob}
below.



\paragraph{Asymptotic equivalence and contiguity.}
\label{sec-asy-eq}

We now define two random graph models that are closely related to
the Norros-Reittu random graph.
In the \emph{generalized random graph model} \cite{BriDeiMar-Lof05},
which we denote by $\GRGnw$, the edge probability of the edge
between vertices $i$ and $j$ is equal to $p_{ij}^{\sss\rm(GRG)}=\frac{w_iw_j}{\ell_n+w_iw_j}.$
In the \emph{random graph with prescribed expected degree} or Chung-Lu random graph
\cite{ChuLu02b,ChuLu03,ChuLu06c,ChuLu06}, which we denote by $\CLnw$,
the edge probabilities are given by $p_{ij}^{\sss({\rm CL})}=\max\{\frac{w_iw_j}{\ell_n},1\}.$
The Chung-Lu model is
sometimes referred to as the random graph with given expected degrees,
as the expected degree of vertex $j$ is close to $w_j$.

By \cite[Examples 3.5 and 3.6]{Jans08a}, in our setting, the graphs
$\NRnw,\GRGnw,$ and $\CLnw$ are \emph{asymptotically equivalent}
(meaning that all events have the same asymptotic probabilities),
so that Theorems \ref{main-clusterdistr-tau>4}--\ref{main-clusterdistr-tau(3,4)}
apply to $\GRGnw$ and $\CLnw$ as well.

\paragraph{The configuration model.}
Given  a {\bf degree sequence}, namely, a sequence of $n$ positive integers
$\bfdit = (d_i)_{i\in [n]}$ with the total degree $\ell_n^{\sss({\rm CM})}=\sum_{i=1}^n d_i$
assumed to be even, the configuration model (CM) on $n$ vertices with
degree sequence $\bfdit$ is constructed as follows:

Start with $n$ vertices and $d_j$ half-edges adjacent to vertex $j$.
Number the half-edges from $1$ to $\ell_n^{\sss({\rm CM})}$ in some arbitrary order.
At each step, two half-edges (not already paired) are chosen uniformly at
random among all the unpaired half-edges and are paired to form a single
edge in the graph. Remove the paired half-edges from the list of unpaired
half-edges. Continue with this procedure until all half-edges are paired.

By varying the degree sequence $\bfdit$, one obtains random graphs with
various degree sequences in a similar way as how varying $\bfwit$ influences
the degree sequence in the $\NRnw$ model studied here.
A first setting which produces a random graph with asymptotic degree sequences
according to some distribution $F$ arises by taking $(d_i)_{i\in[n]}=(D_i)_{i\in[n]}$,
where $(D_i)_{i\in[n]}$ are i.i.d.\ random variables with distribution function $F$.
An alternative choice is to take $(d_i)_{i\in[n]}$ such that the number of
vertices with degree $k$ equals $\lceil nF(k)\rceil - \lceil nF(k-1)\rceil$.

The graph generated in the construction of the
CM is not necessarily \emph{simple}, i.e., it can have self-loops and
multiple edges. However, if
    \eqn{
    \lbeq{conv-nu-CM}
    \nu_n^{\sss({\rm CM})}=\frac{1}{\ell_n^{\sss({\rm CM})}}\sum_{i=1}^n d_i(d_i-1)
    }
converges as $n\rightarrow \infty$
and $d_j=o(\sqrt{n})$ for each $j\in [n]$, then the number of self-loops and multiple
converge in distribution to independent Poisson random variables
(see e.g., \cite{Jans06b} and the references therein). In \cite{MolRee95},
the phase transition of the CM was investigated, and it was shown that
when $\nu_n^{\sss({\rm CM})}\rightarrow \nu^{\sss({\rm CM})}>1$,
and certain conditions on the degrees are satisfied, then a giant component exists,
while if $\nu^{\sss({\rm CM})}\leq 1$, then the largest connected component has size $\op(n)$.
In \cite{JanLuc07}, some of the conditions were removed. Also the
{\it barely supercritical} regime, where $n^{1/3}(\nu_n-1)\rightarrow \infty$,
is investigated. One of the conditions in \cite{JanLuc07} is
that $\sum_{i=1}^n d_i^{4+\eta}=O(n)$ for some $\eta>0$, which,
in the power-law setting, corresponds to $\tau>5$.
Here we write that $f(n)=O(g(n))$ for a non-negative function $g(n)$
when there exists a constant $C>0$ such that $|f(n)|/g(n)\leq C$ for all $n\geq 1$.

In \cite[Remark 2.5]{JanLuc07}, it is conjectured
that this condition is not necessary, and that, in fact, the results should
hold when $\frac 1n \sum_{i=1}^n d_i^{3}$ converges.
Similar results are proved in \cite{KanSei08} under related
conditions. The results in \cite{JanLuc07, KanSei08} suggest that the
barely supercritical regime for the CM is similar to the one for the
ERRG when $\tau>4$. We strengthen
this by conjecturing that Theorems
\ref{main-clusterdistr-tau>4}-\ref{main-clusterdistr-tau(3,4)}
also hold for the CM when $\vep_n=\nu_n-1$ is replaced by
$\vep_n=\nu_n^{\sss({\rm CM})}-1$. After the completion of the first
version of this paper, a
result in this direction was established in \cite{HatMol09}.
Indeed, denote
    \eqn{
    R=\frac1n \sum_{i\in [n]} d_i^3.
    }
Then, under the assumption that the maximal degree of the graph satisfies
$\Delta_n\equiv \max_{i\in [n]} d_i\leq n^{1/3} R^{1/3}/\log{n}$,
and for $|\vep_n|=|\nu_n^{\sss({\rm CM})}-1|\leq \Lambda n^{-1/3} R^{2/3}$,
$|\Cmax|$ is with high probability in between $\omega n^{2/3} R^{-1/3}$ and
$n^{2/3} R^{-1/3}/\omega$ for large $\omega>1$.

When $\tau>4$, $R$ remains uniformly bounded,
while for $\tau\in (3,4)$, $R=\Theta(n^{1-3/(\tau-1)})$ and
$\Delta_n=\Theta(n^{1/(\tau-1)})$. Therefore, for $\tau\in (3,4)$,
the bound on $\Delta_n$ has an extra $1/\log{n}$ too many
to be able to compare it to our results. The high degree vertices
play a crucial role in the scaling limit (as shown
in \cite{BhaHofLee09b}), so that we conjecture that the scaling
limit is affected by this restriction. We also refer
to the discussion in \cite[Section 1.2]{HatMol09} for a discussion
on the relations between our results.  It would be of interest to
investigate the relations further.

It is well know that when \refeq{conv-nu-CM} holds and $\Delta_n=o(\sqrt{n})$,
the CM is asymptotically contiguous to a uniform random graph with the same degree sequence.
Indeed, the number of self-loops and multiple edges converge in distribution
to independent Poisson random variables, which are both equal to zero with positive
probability. Further, the CM conditioned on not having any self-loops is
a uniform random graph with the same degree sequence.
Also the generalized random graph conditioned on its degrees is
also a uniform random graph with that degree sequence (see
e.g., \cite{BriDeiMar-Lof05}). Since the degrees of
the barely supercritical regime in the rank-1 inhomogeneous random graph as
studied here satisfy the conditions in \cite{JanLuc07}, the results there
also apply to our model, whenever $\tau>5$. We leave further details of this argument
to the reader.



\section{Strategy of the proof}
\label{sec-overview}
In this section, we describe the strategy of proof for Theorems
\ref{main-clusterdistr-tau>4}--\ref{main-clusterdistr-tau(3,4)}.
We start by discussing the relevant \emph{first and second moment
methods} in Section \ref{sec-FSMM}, and in Section \ref{sec-reduction},
we reduce the proof to two key propositions.

\subsection{First and second moment methods for cluster sizes}
\label{sec-FSMM}
We denote by
    \eqn{
    \lbeq{Zgeqk-def}
    Z_{\sss \geq k}=\sum_{v\in [n]} \indic{|\cluster(v)|\geq k}
    }
the number of vertices that are contained in connected components
of size at least $k$. Here, we write
$\indicwo{A}$ for the indicator of the event $A$.

The random variable $Z_{\sss \geq k}$ will be used
to prove the asymptotics of $|\Cmax|$. This can be understood by noting
that $|\Cmax|\geq k$ occurs precisely when $Z_{\sss \geq k}\geq k$,
which allows us to prove bounds on $|\Cmax|$ by investigating $Z_{\sss \geq k}$
for appropriately chosen values of $k$. This strategy has been successfully
applied in several related settings, such as percolation on the torus in
general dimension \cite{BorChaKesSpe99} as well as for percolation
on high-dimensional tori \cite{BorChaHofSlaSpe05a,HeyHof06,HofLuc06}.
This is the first time that this methodology is applied to an
\emph{inhomogeneous} setting.

The main aim of this section is to formulate
the necessary bounds on cluster tails and expected cluster size that ensure
the asymptotics in Theorems \ref{main-clusterdistr-tau>4}-\ref{main-clusterdistr-tau(3,4)}.
This will be achieved in Propositions \ref{thm-UBcrit}--\ref{thm-LBcrit} below,
which derive the necessary bounds for the upper and lower bounds on the maximal
cluster size respectively. Throughout the paper, we will use the notation
$(x\wedge y)=\min\{x,y\}, (x\vee y)=\max\{x,y\}$.

\begin{proposition}[An upper bound on the largest critical cluster]
\label{thm-UBcrit}
Fix $\Lambda>0$, and suppose that there exist $\delta>1$ and $a_1=a_1(\Lambda)>0$ such that,
for all $k\geq n^{\delta/(1+\delta)}$ and
for $\Ver$ a uniformly chosen vertex in $[n]$, the bound
    \eqn{
    \lbeq{clustertail-assUB}
    \prob(|\cluster(\Ver)|\geq k)\leq
    a_1\big(k^{-1/\delta}+\big(\vep_n\vee n^{-(\delta-1)/(\delta+1)}\big)^{1/(\delta-1)}\big)
    }
holds, where
    \eqn{
    \lbeq{alpha-restr}
    |\vep_n|\leq
    \Lambda n^{-(\delta-1)/(\delta+1)}.
    }
Then, there exists a $b_1=b_1(\Lambda)>0$ such that, for all $\omega\geq 1$,
    \eqn{
    \prob\big(|\Cmax|\geq \omega n^{\delta/(1+\delta)}\big)\leq \frac{b_1}{\omega}.
    }
\end{proposition}

The bound in \refeq{clustertail-assUB} can be understood as a bound on
the tail of the total progeny of a branching process, where the first term
corresponds to the total progeny being finite and larger than $k$, while the
second term corresponds to the survival probability of the branching process.
This will be made precise in the sequel.

\proof We use the first moment method or Markov inequality, to bound
    \eqn{
    \prob(|\Cmax|\geq k)
    =\prob(Z_{\sss \geq k}\geq k)\leq \frac{1}{k}
    \expec[Z_{\sss \geq k}]=\frac{n}{k}\prob(|\cluster(\Ver)|\geq k),
    }
where $\Ver\in [n]$ is a uniformly chosen vertex.
Thus, we need to bound $\prob(|\cluster(\Ver)|\geq k)$ for an appropriately
chosen $k=k_n$. We use \refeq{clustertail-assUB}, so that
    \eqan{
    \prob(|\Cmax|\geq k)
    &\leq \frac{a_1n}{k}\Big(k^{-1/\delta}+\big(\vep_n\vee n^{-(\delta-1)/(\delta+1)}\big)^{1/(\delta-1)}\Big)\nn\\
    &\leq a_1\big(\omega^{-(1+1/\delta)}+(n^{1/(\delta+1)}|\vep_n|^{1/(\delta-1)}\vee 1)\omega^{-1}\big)\nn\\
    &\leq a_1\big(\omega^{-(1+1/\delta)}+(1+\Lambda^{1/(\delta-1)})\omega^{-1}\big),
    }
when $k=k_n=\omega n^{1/(1+1/\delta)}=\omega n^{\delta/(1+\delta)}$,
and where we have used \refeq{alpha-restr}.
This completes the proof of Proposition \ref{thm-UBcrit}, with $b_1=a_1(2+\Lambda^{1/(\delta-1)})$.
\qed
\vskip0.5cm

\noindent
Proposition \ref{thm-UBcrit} shows that to prove an
upper bound on $|\Cmax|$, it suffices to prove an upper bound
on the cluster tails of a uniformly chosen vertex. In order to prove a
matching lower bound on $|\Cmax|$, we shall use the \emph{second moment
method}, for which we need to give a bound on the variance of
$Z_{\sss \geq k}$. To state the result, we define
    \eqn{
    \lbeq{chi>-def}
    \chi_{\sss \geq k}(\bfpit)=\expec[|\cluster(\Ver)|\indic{|\cluster(\Ver)|\geq k}],
    }
where $\bfpit=(p_{ij})_{1\leq i<j\leq n}$ denote
the edge probabilities of an inhomogeneous random graph, i.e.,
the edge $ij$ is occupied with probability $p_{ij}$ and the
occupation status of different edges are independent.
Then the main variance estimate on $Z_{\sss \geq k}$ is as follows:

\begin{proposition}[A variance estimate for $Z_{\sss \geq k}$]
\label{prop-varZ1}
For any inhomogeneous random graph with edge probabilities
$\bfpit=(p_{ij})_{1\leq i<j\leq n}$, every $n$ and $k\in [n]$,
$\Var(Z_{\sss \geq k})\leq n\chi_{\sss \geq k}(\bfpit).$
\end{proposition}

\proof
We use the fact that
    \eqan{
    {\rm Var}(Z_{\sss \geq k})
    &=\sum_{i,j\in [n]}\big[\prob(|\cluster(i)|\geq k, |\cluster(j)|\geq k)-
    \prob(|\cluster(i)|\geq k)\prob(|\cluster(j)|\geq k)\big].
    }
We split the probability $\prob(|\cluster(i)|\geq k, |\cluster(j)|\geq k)$,
depending on whether $i\conn j$ or not, i.e., we split
    \eqan{
    \prob(|\cluster(i)|\geq k, |\cluster(j)|\geq k)
    &=\prob(|\cluster(i)|\geq k, |\cluster(j)|\geq k, i\conn j)\nn\\
    &\qquad+\prob(|\cluster(i)|\geq k, |\cluster(j)|\geq k, i\nc j).
    }
We can bound
    \eqn{
    \lbeq{disjoccBK}
    \prob(|\cluster(i)|\geq k, |\cluster(j)|\geq k, i\nc j)
    \leq \prob\Big(\{|\cluster(i)|\geq k\}\circ \{|\cluster(j)|\geq k\}\Big),
    }
where, for two increasing events $E$ and $F$, we write $E\circ F$
to denote the event that $E$ and $F$ occur \emph{disjointly}, i.e.,
that there exists a (random) set of edges $K$ such that we can see
that $E$ occurs by only inspecting the edges in $K$ and that
$F$ occurs by only inspecting the edges in $K^c$. Then,
the BK-inequality \cite{BerKes85,Grim99} states that
    \eqn{
    \prob(E\circ F)\leq \prob(E)\prob(F).
    }
Applying this to \refeq{disjoccBK}, we obtain that
    \eqn{
    \prob(|\cluster(i)|\geq k, |\cluster(j)|\geq k, i\nc j)
    \leq \prob(|\cluster(i)|\geq k)\prob(|\cluster(j)|\geq k).
    }
Therefore,
    \eqn{
    {\rm Var}(Z_{\sss \geq k})
    \leq \sum_{i,j\in [n]} \prob(|\cluster(i)|\geq k, |\cluster(j)|\geq k, i\conn j),
    }
and we arrive at the fact that
    \eqan{
    {\rm Var}(Z_{\sss \geq k})
    &\leq \sum_{i,j\in [n]}\prob(|\cluster(i)|\geq k, |\cluster(j)|\geq k, i\conn j)\\
    &= \sum_{i\in [n]}\sum_{j\in [n]}\expec\big[\indic{|\cluster(i)|\geq k}\indic{j\in \cluster(i)}\big]= \sum_{i\in [n]}\expec\Big[\indic{|\cluster(i)|\geq k}\sum_{j\in [n]}\indic{j\in \cluster(i)}\Big]\nn\\
    &=\sum_{i\in [n]}\expec[|\cluster(i)|\indic{|\cluster(i)|\geq k}]=n\expec[|\cluster(\Ver)|\indic{|\cluster(\Ver)|\geq k}]
    =n\chi_{\sss \geq k}(\bfpit).\nn\qed
    }

\begin{proposition}[A lower bound on the largest critical cluster]
\label{thm-LBcrit}
Suppose that there exist $\delta>1$ and $a_2>0$ such that
for all $k\leq n^{\delta/(1+\delta)}$ and
for $\Ver$ a uniformly chosen vertex in $[n]$,
    \eqn{
    \lbeq{clustertail-assLB}
    \prob(|\cluster(\Ver)|\geq k)\geq \frac{a_2}{k^{1/\delta}},
    }
while
    \eqn{
    \lbeq{susc-assLB}
    \expec[|\cluster(\Ver)|]\leq n^{(\delta-1)/(\delta+1)},
    }
then there exists an $b_2>0$ such that, for all $\omega\geq 1$,
    \eqn{
    \prob\big(|\Cmax|\leq \omega^{-1} n^{\delta/(1+\delta)}\big)\leq
    \frac{b_2}{\omega^{2/\delta}}.
    }
\end{proposition}
%

\noindent
\proof We use the second moment method or
Chebychev inequality,
as well as the fact that $|\Cmax|<k$ precisely when $Z_{\sss \geq k}=0$,
to obtain that
    \eqn{
    \lbeq{bdexpecZcrita}
    \prob\big(|\Cmax|<\omega^{-1} n^{\delta/(1+\delta)}\big)
    =\prob\big(Z_{\sss \geq \omega^{-1} n^{\delta/(1+\delta)}}=0\big)
    \leq \frac{\Var(Z_{\sss \geq \omega^{-1} n^{\delta/(1+\delta)}})}
    {\expec[Z_{\sss \geq \omega^{-1} n^{\delta/(1+\delta)}}]^2}.
    }
By \refeq{clustertail-assLB}, we have that
    \eqn{
    \lbeq{bdexpecZcrit}
    \expec[Z_{\sss \geq \omega^{-1} n^{\delta/(1+\delta)}}]
    =n\prob(|\cluster(\Ver)|\geq \omega^{-1} n^{\delta/(1+\delta)})\geq
    \frac{na_2 \omega^{1/\delta}}{n^{1/(1+\delta)}}= a_2\omega^{1/\delta} n^{\delta/(\delta+1)}.
    }
Also, by Proposition \ref{prop-varZ1}, with $k=k_n=\omega^{-1} n^{\delta/(\delta+1)}$,
    \eqn{
    \lbeq{bdvarZcrit}
    \Var(Z_{\sss \geq \omega^{-1} n^{\delta/(\delta+1)}})\leq
    n\chi_{\sss \geq \omega^{-1} n^{\delta/(\delta+1)}}(\bfpit)
    \leq n^{1+(\delta-1)/(\delta+1)}=n^{2\delta/(\delta+1)}.
    }
Substituting \refeq{bdexpecZcrita}--\refeq{bdvarZcrit}, we obtain,
for $n$ sufficiently large,
    \eqn{
    \prob_{1}\big(|\Cmax|<\omega^{-1} n^{\delta/(1+\delta)}\big)
    \leq \frac{n^{2\delta/(\delta+1)}}{a_2^2\omega^{2/\delta} n^{2\delta/(\delta+1)}}
    =\frac{1}{a_2^2\omega^{2/\delta}}.}
This completes the proof of Proposition \ref{thm-LBcrit}.
\qed

\subsection{Reduction of the proof to two key propositions}
\label{sec-reduction}
In this section, we state two key proposition and use it to complete the proof of
Theorems \ref{main-clusterdistr-tau>4}--\ref{main-clusterdistr-tau(3,4)}.
We denote
    \eqn{
    \lbeq{delta-def}
    \delta=(\tau\wedge 4)-2.
    }
Our main technical results are formulated in the following two propositions:

\begin{proposition}[An upper bound on the cluster tail]
\label{prop-main}
Fix $\NRnw$ with $\bfwit=(w_j)_{j\in [n]}$ as in \refeq{choicewi}, and assume that
the distribution function $F$ in \refeq{choicewi} satisfies $\nu=1$.
Assume that \refeq{F-bound-tau>4} holds for some $\tau>4$, or that
\refeq{F-bound-tau(3,4)b} holds for some $\tau\in (3,4)$, and fix
$\vep_n$ such that $|\vep_n|\leq \Lambda n^{-(\delta-1)/(\delta+1)}$
for some $\Lambda>0$. Let $\tilde{\bfwit}$ be defined as in
\refeq{tilde-wi-def}. Then, for all $k\geq 1$ and
for $\Ver$ a uniformly chosen vertex in $[n]$:\\
(a) There exists a constant $a_1>0$ such that
    \eqn{
    \lbeq{clustertail-UB}
    \prob(|\cluster(\Ver)|\geq k)\leq
    a_1\big(k^{-1/\delta}+\big(\vep_n\vee n^{-(\delta-1)/(\delta+1)}\big)^{1/(\delta-1)}\big).
    }
(b) There exists a constant $a_2>0$ such that
    \eqn{
    \lbeq{clustertail-LB}
    \prob(|\cluster(\Ver)|\geq k)\geq \frac{a_2}{k^{1/\delta}}.
    }
\end{proposition}

\begin{proposition}[An upper bound on the expected cluster size]
\label{thm-UBclust-size}
Fix $\Lambda\geq 1$ sufficiently large, and
let $\vep_n\leq -\Lambda n^{-(\delta-1)/(\delta+1)}$.
Then, for $\Ver$ a uniformly chosen vertex in $[n]$,
    \eqn{
    \lbeq{susc-assLB-proved}
    \expec[|\cluster(\Ver)|]\leq n^{(\delta-1)/(\delta+1)}.
    }
\end{proposition}

\noindent
Now we are ready to prove Theorems \ref{main-clusterdistr-tau>4}--\ref{main-clusterdistr-tau(3,4)}:

\paragraph{Upper bounds in Theorems \ref{main-clusterdistr-tau>4}--\ref{main-clusterdistr-tau(3,4)}.}
The upper bounds follow immediately from Propositions \ref{thm-UBcrit} and \ref{prop-main}(a),
when we recall the definition of $\delta$ in
\refeq{delta-def}, so that \refeq{alpha-restr} is the same as
$|\vep_n|\leq \Lambda n^{-1/3}$ for $\tau>4$, as assumed in Theorem \ref{main-clusterdistr-tau>4},
and $|\vep_n|\leq \Lambda n^{-(\tau-3)/(\tau-2)}$ for $\tau\in(3,4)$,
as assumed in Theorem \ref{main-clusterdistr-tau(3,4)}.

\paragraph{Lower bounds in Theorems \ref{main-clusterdistr-tau>4}--\ref{main-clusterdistr-tau(3,4)}.}
For the lower bounds, we note that there is an obvious monotonicity in
the weights, so that the cluster for $\tilde{w}_i$ as in \refeq{tilde-wi-def}
with $\tilde\vep_n=-\Lambda n^{-(\delta-1)/(\delta+1)}$ is stochastically smaller than the one for
$\vep_n$ with $|\vep_n|\leq \Lambda n^{-(\delta-1)/(\delta+1)}$.
Then, we make use of Proposition \ref{thm-LBcrit}, and check that its assumptions
are satisfied due to Propositions \ref{prop-main}(b) and \ref{thm-UBclust-size}.
\qed


\section{Preliminaries}
\label{sec-prel}
In this section, we derive preliminary results needed in the proofs of
Propositions \ref{prop-main} and \ref{thm-UBclust-size}. We start in Section
\ref{sec-moments-w} by analyzing sums of functions of the vertex weights,
and in Section \ref{sec-MPBP} we describe
a beautiful connection between branching processes and clusters
in the Norros-Reittu model originating in \cite{NorRei06}.

\subsection{The weight $W_n$ of a uniformly chosen vertex}
\label{sec-moments-w}
In this section, we investigate the weight of a uniformly chosen vertex in $[n]$, which
we denote by $W_n$. For this, we first note that
    \eqn{
    \lbeq{Finvid}
    [1-F]^{-1}(1-u)=F^{-1}(u)=\inf\{x: F(x)\geq u\},
    }
which, in particular, implies that $W=[1-F]^{-1}(U)$ has
distribution function $F$ when $U$ is uniform on $(0,1)$.
Further, $W_n$ is a random variable with distribution
function $F_n$ given by
    \eqan{
    \lbeq{Fn-def}
    F_n(x)&=\prob(W_n\leq x)=\frac 1n\sum_{j=1}^n \indic{w_j\leq x}=\frac 1n\sum_{j=1}^n \indic{[1-F]^{-1}(\frac{j}{n})\leq x}=\frac 1n\sum_{i=0}^{n-1} \indic{[1-F]^{-1}(1-\frac{i}{n})\leq x}\nn\\
    &=\frac 1n\sum_{i=0}^{n-1}\indic{F^{-1}(\frac{i}{n}) \leq x}
    =\frac 1n\sum_{i=0}^{n-1} \indic{\frac{i}{n}\leq F(x)}
    =\frac 1n\big(\big\lfloor n F(x)\big\rfloor+1\big)\wedge 1,
    }
where we write $j=n-i$ in the fourth equality and use \refeq{Finvid} in the fifth equality.
Note that $F_n(x)\geq F(x)$, which shows that $W_n$ is stochastically dominated
by $W$, so that, in particular, for \emph{increasing}
functions $x\mapsto h(x)$,
    \eqn{
    \lbeq{sum-h-increasing}
    \frac{1}{n}\sum_{j=1}^n h(w_j)\leq \expec[h(W)].
    }
In the sequel, we shall repeatedly bound
expectations of functions of $W_n$ using the following lemma:

\begin{lemma}[Expectations of $W_n$]
\label{lem-conv-sums-w}
Let $W$ have distribution function $F$ and let $W_n$ have distribution function $F_n$ in
\refeq{Fn-def}. Let $h \colon [0,\infty)\to {\mathbb C}$ be a differentiable
function with $h(0)=0$ such that $|h'(x)|[1-F(x)]$ is integrable on $[0,\infty)$.
Then,
    \eqn{
    \big|\expec[h(W_n)]-\expec[h(W)]\big|
    \leq \int_{w_1}^{\infty} |h'(x)|[1-F(x)]dx
    +\frac{1}{n}\int_0^{w_1}|h'(x)|dx.
    }
\end{lemma}

\proof We write, using the fact that $h(0)=0$ and that $|h'(x)|[1-F(x)]$ is integrable,
for any $B>0$,
    \eqan{
    \lbeq{expec(W)}
    \expec[h(W)\indic{W\leq B}]
    &\expec\Big[\int_{0}^{\infty} h'(x)\indic{x<W\leq B}dx\Big]
    =\int_{0}^{B} h'(x)[F(B)-F(x)]dx.
    }
Because of this representation, we have that
    \eqn{
    \expec[h(W_n)]-\expec[h(W)]
    =\int_{0}^{\infty} h'(x)[F(x)-F_n(x)]dx.
    }
Now, $F_n(w_1)=1$ by construction (recall \refeq{Fn-def}), so that
    \eqn{
    \Big|\expec[h(W_n)]-\expec[h(W)]\Big|
    \leq \int_{w_1}^{\infty} |h'(x)|[1-F(x)]dx
    +\int_0^{w_1}|h'(x)|[F_n(x)-F(x)]dx.
    }
We finally use the fact that $0\leq F_n(x)-F(x)\leq 1/n$ to arrive at the claim.
\qed

\begin{corollary}[Bounds on characteristic function and mean degrees]
\label{cor-varphi-n}
Let $W$ and $W_n$ have distribution functions $F$ and $F_n$, respectively,
and assume that \refeq{F-bound-tau>4} holds
for some $\tau>3$.\\
(a) Let
    \eqn{
    \lbeq{varphi-n-def}
    \phi_n(t)=\frac{\expec[W_n\mathrm{e}^{(1+\vep_n)W_n(\mathrm{e}^{\i t}-1)}]}{\expec[W_n]},
    \qquad
    \phi(t)=\frac{\expec[W\mathrm{e}^{W(\mathrm{e}^{\i t}-1)}]}{\expec[W]}.
    }
Then,
    \eqn{
    |\phi_n(t)-\phi(t)|\leq cn^{-(\tau-2)/(\tau-1)}+c|t|(n^{-(\tau-3)/(\tau-1)}+|\vep_n|).
    }
(b) With $\nu$ as in \refeq{nu-def}, $\nu_n$ as in\refeq{nu-n-def} and with $\tilde\nu_n=(1+\vep_n)\nu_n$,
    \eqn{
    \lbeq{nu-n-conv}
    |\tilde\nu_n-\nu|\leq c(|\vep_n|+n^{-(\tau-3)/(\tau-1)}).
    }
\end{corollary}
We remark that if \refeq{F-bound-tau>4} holds
for some $\tau>4$ or \refeq{F-bound-tau(3,4)b} holds
for some $\tau\in(3,4)$, then also \refeq{F-bound-tau>4} holds
for that $\tau>3$. This explains the assumption that $\tau>3$ in
Corollary \ref{cor-varphi-n}.

\proof (a) We first take $\vep_n=0$, and split
    \eqn{
    \phi_n(t)-\phi(t)
    =\frac{\phi(t)}{\expec[W]}\Big(\frac{1}{\expec[W_n]}-\frac{1}{\expec[W]}\big)
    +\frac{1}{\expec[W_n]}\Big(\expec[W_n\mathrm{e}^{W_n(\e^{\i t}-1)}]-\expec[W\mathrm{e}^{W(\e^{\i t}-1)}]\Big).
    }
Lemma \ref{lem-conv-sums-w} applied to $h(x)=x$ yields
$\big|\expec[W_n]-\expec[W]\big| \leq cn^{-(\tau-2)/(\tau-1)}.$
To apply Lemma \ref{lem-conv-sums-w} to $h(x)=x\mathrm{e}^{x(\e^{\i t}-1)},$ we compute
    \eqn{
    |h'(x)|=|\mathrm{e}^{x(\e^{\i t}-1)}+i(\e^{\i t}-1)x\mathrm{e}^{it(x-1)}|
    \leq 1+x|t|.
    }
Therefore, also using the fact that $w_1=\Theta(n^{1/(\tau-1)})$ by \refeq{F-bound-tau>4},
    \eqan{
    \Big|\expec[W_n\mathrm{e}^{W_n(\e^{\i t}-1)}]-\expec[W\mathrm{e}^{W(\e{\i t}-1)}]\Big|
    &\leq \int_{w_1}^{\infty} (1+x|t|)[1-F(x)]dx+\frac{1}{n}\int_0^{w_1}(1+x|t|)dx\nn\\
    &\leq cn^{-(\tau-2)/(\tau-1)}+c|t|n^{-(\tau-3)/(\tau-1)}.
    }
Together, these two estimates prove the claim for $\vep_n=0$.
For $\vep_n\neq 0$, we use the fact that
    \eqan{
    \frac{\expec[W_n\mathrm{e}^{(1+\vep_n)W_n(\e^{\i t}-1)}]}{\expec[W_n]}
    -\frac{\expec[W_n\mathrm{e}^{W_n(\e^{\i t}-1)}]}{\expec[W_n]}
    &=\frac{\expec\big[W_n\mathrm{e}^{W_n(\e^{\i t}-1)}(\mathrm{e}^{\vep_n W_n(\e^{\i t}-1)}-1)\big]}{\expec[W_n]}\nn\\
    &=O(|\vep_n||t|\expec[W_n^2]/\expec[W_n])=O(|\vep_n||t|).
    }
(b) The proof of \refeq{nu-n-conv} is similar.
\qed

\begin{lemma}[Bounds on moments of $W_n$]
\label{lem-MomtailF}
Let $W_n$ have distribution function $F_n$ in \refeq{Fn-def}.
\begin{itemize}
\item[(i)] Assume that
\refeq{F-bound-tau>4} holds for some $\tau>3$, and let $a<\tau-1$. Then,
for $x$ sufficiently large, there exists a $C=C(a,\tau)$ such that,
uniformly in $n$,
    \eqn{
    \expec[W^a_n\indic{W_n\geq x}]\leq Cx^{a+1-\tau}.
    }
\item[(ii)] Assume that
\refeq{F-bound-tau(3,4)b} holds for some $\tau>3$, and let $a>\tau-1$.
Then, there exist $C_1=C_1(a,\tau)$ and $C_2=C_2(a,\tau)$ such that, uniformly in $n$,
    \eqn{
    C_1\big(x\wedge n^{1/(\tau-1)}\big)^{a+1-\tau}\leq \expec[W^a_n\indic{W_n\leq x}]\leq C_2x^{a+1-\tau}.
    }
\end{itemize}
\end{lemma}

\proof (i) When $a<\tau-1$, the expectation is finite. We rewrite, using \refeq{expec(W)},
    \eqan{
    \expec[W_n^a\indic{W_n\geq x}]
    &=x^a[1-F_n](x) + a\int_{x}^{\infty} v^{a-1}[1-F_n(v)]dv.
    }
Now, $1-F_n(x)\leq 1-F(x)$, so that
we may replace the $F_n$ by $F$ in an upper bound. When \refeq{F-bound-tau>4} holds
for some $\tau>3$, we can further bound this as
    \eqn{
    \expec[W_n^a\indic{W\geq x}]\leq c_{\sss F} x^{a+1-\tau}+c_{\sss F} a\int_{x}^{\infty} w^{a-\tau}dw
    =O(x^{a+1-\tau}).
    }
(ii) We again use \refeq{expec(W)} and the bound in \refeq{F-bound-tau(3,4)b} to rewrite
    \[
    \expec[W^a_n\indic{W_n\leq x}]
    =a\int_{0}^{x} v^{a-1}[F_n(x)-F_n(v)]dv\leq a\int_{0}^{x} v^{a-1}[1-F(v)]dv
    \leq c_{\sss F} a\int_{0}^{x} v^{a-\tau}dv=O(x^{a+1-\tau}).
    \]
For the lower bound, we first assume that $x\leq n^{1/(\tau-1)}$ and use the fact that
    \eqan{
    \expec[W_n^a\indic{W_n\leq x}]&=a\int_{0}^{x} v^{a-1}[F_n(x)-F_n(v)]dv
    \geq a\int_{0}^{\vep x} v^{a-1}[F_n(x)-F_n(v)]dv\nn\\
    &=a\int_{0}^{\vep x} v^{a-1}[1-F_n(v)]dv-a^{-1}[1-F_n(x)]\int_{0}^{\vep x} v^{a-1}dv\nn\\
    &\geq a\int_{0}^{\vep x} v^{a-1}[1-F(v)]dv-\frac{a}{n}\int_{0}^{\vep x}v^{a-1} dv
    -a[1-F(x)]\int_{0}^{\vep x} v^{a-1}dv\nn\\
    &\geq C(\vep x)^{a+1-\tau}-C(\vep x)^a/n-Cx^{-(\tau-1)}(\vep x)^a\nn\\
    &=C \vep^{a+1-\tau} x^{a+1-\tau}\Big(1-\vep^{\tau-1}\big(x/n^{1/(\tau-1)}\big)^{\tau-1}-\vep^{\tau-1}\Big)
    \geq C_2 x^{a+1-\tau},
    }
when we take $\vep\in (0,1)$ sufficiently small, and we use the fact that $x\leq n^{1/(\tau-1)}$. When
$x\geq n^{1/(\tau-1)}$, we can use the fact that $w_1=\Theta(n^{1/(\tau-1)})$,
so that
    \eqn{
    \expec[W_n^a\indic{W_n\leq x}]\geq w_1^a/n\geq (c n^{1/(\tau-1)})^a/n=
    C_2 \big(n^{1/(\tau-1)}\big)^{a+1-\tau}.
    }
\qed

\subsection{Connection to mixed Poisson branching processes}
\label{sec-MPBP}
In this section, we discuss the relation between our Poisson random graph
and mixed Poisson branching processes due to Norros and Reittu \cite{NorRei06}.


\paragraph{Stochastic domination of clusters
by a branching process.}
We shall dominate the cluster of a vertex in the
Norros-Reittu model by the total progeny of an appropriate
branching process. In order to describe this relation,
we consider the cluster exploration
of a uniformly chosen vertex $\Ver\in[n]$. For this, we
define the \emph{mark distribution} to be the random variable
$M$ with distribution
    \eqn{
    \lbeq{dist-mark}
    \prob(M=m)= w_m/\ell_n, \qquad
    m\in [n].
    }
We define $S_0=1$, and, recursively, for $i\geq 1$,
    \eqn{
    \lbeq{Si-rec}
    S_i=S_{i-1}+X_i-1,
    }
where $(X_i)_{i\geq 1}$ is a sequence of independent random variables, where
$X_i$ has a mixed Poisson distribution with random parameter $w_{\sss M_i}$,
and where $M_1$ is uniformly chosen in $[n]$, while
$(M_i)_{i\geq 2}$ are i.i.d.\ random marks with distribution
\refeq{dist-mark}. Let
    \eqn{
    \lbeq{TP-BP}
    T^{\sss(2)}=\inf\{t\geq 1\colon S_t=0\}
    }
denote the first hitting time of 0 of $(S_i)_{i\geq 0}$.
By exploring a branching process tree, we see that $T^{\sss(2)}$ has the same distribution
as the total progeny of a so-called \emph{two-stage mixed Poisson branching process},
in which the root has $X_1\sim \Poi(w_{\sss \Ver})$ children where $\Ver$ is
chosen uniformly in $[n]$, and all other individuals have
offspring distribution given by $\Poi(w_{\sss M})$. In the sequel,
we shall use the notation $T^{\sss(2)}$ for the total progeny of
a two-stage mixed Poisson branching process with offspring distribution
$X_i\sim \Poi(w_{\sss M_i})$ and the root has offspring distribution
$X_1\sim \Poi(w_{\sss \Ver})$. We shall use the notation $T$ for the total
progeny of a mixed Poisson branching process where every individual,
including the root, has an i.i.d.\ offspring distribution $X_i\sim \Poi(w_{\sss M_i})$,
where $(M_i)_{i\geq 1}$ is an i.i.d.\ sequence of marks with distribution described in
\refeq{dist-mark}.

Clearly, $w_{\Ver}$ has distribution $W_n$
defined in \refeq{Fn-def}, while
    \eqn{
    \lbeq{ml_M}
    \prob(w_{\sss M}\leq x)=\sum_{m=1}^n \indic{w_m\leq x}\prob(M=m)
    =\frac{1}{\ell_n}\sum_{m=1}^n w_m\indic{w_m\leq x}
    =\prob(W_n^*\leq x),
    }
where $W_n^*$ is the \emph{size-biased distribution} of $W_n$ (recall \refeq{size-biased-distr}).
In order to define the cluster exploration in $\NRnw$, we define $\tildeX_1=X_1$ and, for $i\geq 2$,
    \eqn{
    \lbeq{tildeX-def}
    \tildeX_i= X_i \indic{M_i\not\in \{M_1, \ldots, M_{i-1}\}},
    }
and define $\tildeS_i$ and $\tildeT$ as in \refeq{Si-rec} and
\refeq{TP-BP}, where $X_i$ is replaced by $\tildeX_i$.
The definition in \refeq{tildeX-def} can be thought of
as a \emph{thinning} of the branching process. The mark $M_i$
corresponds to the \emph{vertex label} corresponding to the
$i^{\rm th}$ individual encountered in the exploration process.
If we have already explored this vertex, then we should ignore it
and all of its subsequent offspring, while if it is a new vertex,
then we should keep it and explore its offspring. Thus,
the definition in \refeq{tildeX-def} can be thought of ensuring
that we only explore the offspring of a vertex once.
We think of $\tildeX_i$ as the \emph{potential} new vertices
of the cluster $\cluster(\Ver)$ neighboring the $i^{\rm th}$ explored vertex.
A potential vertex turns into a real new element of $\cluster(\Ver)$
when its mark or vertex label is one that we have not yet seen.

We will now make the connection between the thinned marked
mixed Poisson branching process and the cluster exploration precise;
this relation (in a slightly different context) was first
proved in \cite[Proposition 3.1]{NorRei06}. Since
the proof of \cite[Proposition 3.1]{NorRei06} is in terms
of the number of vertices at distance $k$ of the root, we reprove this result
here in the setting of the cluster exploration:

\begin{proposition}[Clusters as thinned marked branching processes]
\label{prop-N NR kop}
The cluster $\cluster(\Ver)$ of a uniformly chosen vertex in $[n]$
is equal in distribution to $\{M_1, \ldots, M_{\sss\tildeT}\}$,
i.e., the marks encountered in the thinned marked mixed Poisson branching process
up to the end of the exploration $\tildeT$.
\end{proposition}

\proof We start by proving that the direct neighbors of the root agree in both
constructions. We note that $\tildeX_1= X_1$, which has a mixed Poisson distribution
with mixing distribution $w_{\sss \Ver}$, which has the same distribution as $W_n$.

Conditionally on $\Ver=l$, $X_1$ has a $\Poi(w_l)$ distribution.
These $\Poi(w_l)$ offspring receive i.i.d.\ marks.
As a result, the random vector $(\tildeX_{1,j})_{j\in [n]}$, where $\tildeX_{1,j}$ is
the number of offspring of the root that receive mark $j$, is a vector of \emph{independent}
Poisson random variables with parameters $w_lw_k/\ell_n$. Due to the thinning,
a mark occurs precisely when $\tildeX_{1,j}\geq 1$, and these events are independent.
Therefore, the mark $j$ occurs, independently for all $j\in [n]$, with probability
$1-\e^{w_jw_k/\ell_n}=p^{\sss{\rm (NR)}}_{jk}$, which proves that the set of neighbors
of the root is equal in distribution to the marks found in our branching process.

Next, we look at the number of new vertices of $\cluster(\Ver)$ neighboring
the $i^{\rm th}$ explored potential vertex. First,
conditionally on $M_i=l$, and assume that $l\not\in \{M_1, \ldots, M_{i-1}\}$,
so that $\tildeX_i= X_i$. Conditionally on $M_i=j$ such that
$j\not\in \{M_1, \ldots, M_{i-1}\}$, $\tildeX_i$ has a $\Poi(w_J)$ distribution.
Each of these $\Poi(w_j)$ potential new vertices in $\cluster(\Ver)$
receives an i.i.d.\ mark. Let $\tildeX_{i,j}$ denote
the number of potential new vertices of $\cluster(\Ver)$
that receive mark $j$. Then, $(\tildeX_{i,j})_{j\in [n]}$
again is a vector of \emph{independent}
Poisson random variables with parameters $w_lw_k/\ell_n$.
Due to the thinning, a mark occurs precisely when $\tildeX_{i,j}\geq 1$,
and these events are independent. In particular,
for each $j\not\in \{M_1, \ldots, M_{i}\}$, the probability that
the mark $j$ occurs equals $1-\e^{w_jw_k/\ell_n}=p^{\sss{\rm (NR)}}_{jk}$,
as required.
\qed

Proposition \ref{prop-N NR kop} implies that, for all $k\geq 1$,
    \eqn{
    \lbeq{SD-cluster}
    \prob(|\cluster(\Ver)|\geq k)\leq \prob(\tildeT\geq k)\leq \prob(T^{\sss(2)}\geq k).
    }
Interestingly, when the weights equal $w_i=n\log{(1-\lambda/n)}$ for each $i\in [n]$,
for which the graph is an Erd\H{o}s-R\'enyi random graph with edge probability
$\lambda/n$, the above implies that $|\cluster(\Ver)|$ is stochastically dominated
by the total progeny of a Poisson branching process with parameter $n\log{(1-\lambda/n)}$.
Often, the stochastic domination is by a branching process with binomial
$n-1$ and $p=\lambda/n$ offspring distribution instead.

\paragraph{Otter-Dwass formula for the branching process total progeny.}
Our proofs make crucial use of
the {\em Otter-Dwass formula}, which describes the distribution of the
total progeny of a branching process (see \cite{Dwas69} for the special case
when the branching process starts with a single individual and
\cite{Otte49} for the more general case, and \cite{HofKea07} for a simple
proof based on induction).

\begin{lemma}[Otter-Dwass formula]
\label{lem-OD}
Let $(X_i)_{i\geq 1}$ be i.i.d.\ random variables. Let $\prob_m$
denote the Galton-Watson process measure with offspring distribution $X_1$
started from $m$ initial individuals, and denote
its total progeny by $T$. Then, for all $k,m\in \N$,
    \eqn{
    \lbeq{OD-form}
    \prob_m(T=k) =\frac{m}{k} \prob (\sum_{i=1}^k X_i = k-m).
    }
\end{lemma}

\paragraph{The survival probability of near-critical mixed Poisson branching processes.}
We shall also need bounds on the survival probability of near-critical
mixed Poisson branching processes.

\begin{lemma}[Survival probability of near-critical mixed Poisson branching processes.]
\label{lem-surv-prob}
Let $\rho_n$ be the survival probability of a branching process with
a $\Poi(W_n^*)$ offspring distribution. Assume that $\vep_n=\expec[W_n^*]-1=\nu_n-1\geq 0$ and
$\vep_n=o(1)$. When \refeq{F-bound-tau>4} holds for some $\tau>4$, there exists a constant
$c>0$ such that
    \eqn{
    \lbeq{rho-n-bd-tau>4}
    \rho_n\leq c \vep_n.
    }
When \refeq{F-bound-tau(3,4)b} holds for some $\tau\in (3,4)$,
    \eqn{
    \lbeq{rho-n-bd-tau(3,4)}
    \rho_n\leq c \big(\vep_n^{1/(\tau-3)}\vee n^{-1/(\tau-1)}\big).
    }
\end{lemma}

\proof By conditioning on the first generation, the survival
probability $\rho$ of a branching process with
offspring distribution $X$ satisfies
    \eqn{
    1-\rho=\expec[(1-\rho)^{X}].
    }
In our case, $X\sim \Poi(W_n^*)$, for which $\expec[s^{\Poi(W_n^*)}]
=\expec[\e^{W_n^*(s-1)}]$, so that
    \eqn{
    1-\rho_n=\expec[\mathrm{e}^{-\rho_n W_n^*}].
    }
We use the fact that $\mathrm{e}^{-x}\geq 1-x$ when $x\geq 1/2$
and $\mathrm{e}^{-x}\geq 1-x+x^2/4$ when $x\leq 1/2$, to arrive at
    \eqn{
    1-\rho_n
    \geq 1-\rho_n\expec[W_n^*] +\frac{\rho_n^2}{4}\expec\big[(W_n^*)^2\indic{\rho_n W_n^*\leq 1/2}\big].
    }
Rearranging terms, dividing by $\rho_n$ and using that
$\expec[W_n^*]=\nu_n$, we obtain
    \eqn{
    \rho_n \expec\big[(W_n^*)^2\indic{\rho_n W_n^*\leq 1/2}\big]\leq 4(\nu_n-1)=4\vep_n.
    }
When \refeq{F-bound-tau>4} holds for some $\tau>4$,
    \eqn{
    \expec\big[(W_n^*)^2\indic{\rho_n W_n^*\leq 1/2}\big]
    =\expec[(W^*)^2](1+o(1))=\frac{\expec[W^3]}{\expec[W]}(1+o(1)),
    }
so that $\rho_n\leq c\vep_n$ for some constant $c>0$. When, on the other hand,
\refeq{F-bound-tau(3,4)b} holds for some $\tau\in (3,4)$, then when
$\rho_n\leq 2n^{-1/(\tau-1)}$ the claimed bound holds.
When, instead, $\rho_n\geq 2n^{-1/(\tau-1)}$,
we may apply Lemma \ref{lem-MomtailF}(ii) to obtain
    \eqn{
    \expec\big[(W_n^*)^2\indic{\rho_n W_n^*\leq 1/2}\big]
    =\frac{\expec\big[W_n^3\indic{\rho_n W_n\leq 1/2}\big]}{\expec[W_n]}
    \geq c\rho_n^{\tau-4},
    }
so that $c\rho_n^{\tau-3}\leq 4\vep_n.$ Combining these two bounds
proves \refeq{rho-n-bd-tau(3,4)}.
\qed

\section{An upper bound on the cluster tail}
\label{sec-UBcluster_tail}
In this section, we shall prove the upper bound on the tail probabilities of
critical clusters in Proposition \ref{prop-main}(a).


\paragraph{Dominating the two-stage branching process by an ordinary branching process.}
We rely on \refeq{SD-cluster}. Unfortunately, the Otter-Dwass formula (Lemma \ref{lem-OD})
is not directly valid for $T^{\sss(2)}$, and we first establish that,
for every $k\geq 0$,
    \eqn{
    \lbeq{T2-SD-T}
    \prob(T^{\sss(2)}\geq k)\leq \prob(T\geq k).
    }
The bound in \refeq{T2-SD-T} is equivalent to
the fact that $T^{\sss(2)}\SD T$, where $X\SD Y$ means that $X$
is stochastically smaller than $Y$. Since the distributions of
$T^{\sss(2)}$ and $T$ agree except for the offspring of the root,
where $T^{\sss(2)}$ has offspring distribution $\Poi(W_n)$, whereas $T$
has offspring distribution $\Poi(W_n^*)$,
this follows when $\Poi(W_n)\SD \Poi(W_n^*)$,
For two mixed Poisson random variables
$X,Y$ with mixing random variables $W_{\sX}$ and $W_{\sY}$,
respectively, $X\SD Y$ follows when $W_{\sX}\SD W_{\sY}$.
%
%
The proof of \refeq{T2-SD-T} is completed by noting that, for any non-negative
random variable $W$, and for $W^*$ its size-biased version, we have $W\SD W^*$.

\paragraph{The total progeny of our mixed Poisson branching process.}
By \refeq{SD-cluster} and \refeq{T2-SD-T},
    \eqan{
    \lbeq{OD-form-applied}
    \prob(|\cluster(\Ver)|\geq k)&\leq \prob(T\geq k)
    =\prob(T=\infty)+\prob(k\leq T<\infty
    =\prob(T=\infty)+\sum_{l=k}^{\infty} \prob(T=l)\nn\\
    &=\prob(T=\infty)+\sum_{l=k}^{\infty} \frac{1}{l} \prob (\sum_{i=1}^l X_i = l-1),
    }
where the last formula follows from Lemma \ref{lem-OD} for $m=1$, and where $(X_i)_{i=1}^{\infty}$
is an i.i.d.\ sequence with a $\Poi(W_n^*)$ distribution. In the following proposition, we
shall investigate $\prob (\sum_{i=1}^l X_i = l-1)$:

\begin{proposition}[Upper bound on probability mass function of $\sum_{i=1}^l X_i$]
\label{prop-sumsXi}
Let $(X_i)_{i=1}^{\infty}$ be an i.i.d.\ sequence with a mixed Poisson distribution
with mixing random variable $\tilde W_n^*=(1+\vep_n)W_n^*$, where $W_n^*$ is defined
in \refeq{ml_M}.
Under the assumptions of Proposition \ref{prop-main}(a), there exists an $\tilde a_1>0$ such that
for all $l\geq n^{\delta/(1+\delta)}$ and for $n$ sufficiently large,
    \eqn{
    \lbeq{sums-Xi-bd}
    \prob \big(\sum_{i=1}^l X_i = l-1\big)\leq \tilde a_1\Big(l^{-1/\delta}+\big(n^{(\tau-4)/2(\tau-1)}\wedge1\big)l^{-1/2}\Big),
    }
where $\delta>0$ is defined in \refeq{delta-def}.
\end{proposition}

\proof We rewrite, using the Fourier inversion theorem, and
recalling $\phi_n(t)$ in \refeq{varphi-n-def}, which we can identify as
$\phi_n(t)=\expec[\mathrm{e}^{\i t X_1}]=\expec[\mathrm{e}^{(\mathrm{e}^{\i t}-1)\tilde W_n^*}]$,
    \eqn{
    \lbeq{FIT-applied}
    \prob \big(\sum_{i=1}^l X_i = l-1\big)=\int_{[-\pi,\pi]} \mathrm{e}^{-\i(l-1)t} \phi_n(t)^{l}\frac{dt}{2\pi},
    }
so that
    \eqn{
    \prob \big(\sum_{i=1}^l X_i = l-1\big)\leq \int_{[-\pi,\pi]} |\phi_n(t)|^{l}\frac{dt}{2\pi}.
    }
By dominated convergence and the weak convergence of $\tilde W_n^*$ to $W^*$, for every $t\in [-\pi,\pi]$,
    \eqn{
    \lbeq{phi-n-conv}
    \lim_{n\rightarrow \infty} \phi_n(t)=\phi(t)=\expec[\mathrm{e}^{W^*(\mathrm{e}^{\i t}-1)}].
    }
Since, further,
    \eqn{
    |\phi_n'(t)|=\big|\expec\big[\tilde W_n^*\mathrm{e}^{\i t}\mathrm{e}^{(\mathrm{e}^{\i t}-1)\tilde W_n^*}\big]\big|
    \leq \expec[\tilde W_n^*]=(1+\vep_n)\nu_n=1+o(1),
    }
which is uniformly bounded, the convergence in \refeq{phi-n-conv} is uniform for
all $t\in[-\pi,\pi]$. Further, for every $\eta>0$, there exists $\vep>0$
such that $|\phi(t)|<1-2\vep$ for all $|t|>\eta$, since
our mixed Poisson random variable is not degenerate at 0.
Therefore, uniformly for sufficiently large $n$,
for every $\eta>0$, there exists $\vep>0$
such that $|\phi_n(t)|<1-\vep$ for all $|t|>\eta$. Thus,
    \eqn{
    \int_{[-\pi,\pi]} |\phi_n(t)|^{l}\frac{dt}{2\pi}
    \leq (1-\vep)^l +\int_{[-\eta,\eta]} |\phi_n(t)|^{l}\frac{dt}{2\pi}.
    }
We start by deriving the bound when $\tau>4$, by bounding
    \eqn{
    |\phi_n(t)|\leq \expec[\mathrm{e}^{-\tilde W_n^*[1-\cos(t)]}].
    }
Now using the fact that, uniformly for $t\in [-\pi,\pi]$, there exists an $a>0$ such that
    \eqn{
    \lbeq{1-cos-bd}
    1-\cos(t)\geq at^2,
    }
and, for $x\leq 1$, the bound $\mathrm{e}^{-x}\leq 1-x/2$, we arrive at
    \eqan{
    |\phi_n(t)| &\leq \expec[\mathrm{e}^{-a\tilde W_n^*t^2}]
    \leq\expec\big[(1-a \tilde W_n^*t^2/2)\indic{a\tilde W_n^*t^2\leq 1}\big]
    +\expec\big[\indic{a\tilde W_n^*t^2>1}\big]\nn\\
    &=1-a t^2 \expec[\tilde W_n^*]
    +\expec\big[\indic{a\tilde W_n^* t^2>1}(1+a\tilde W_n^* t^2/2)\big].
    }
Further bounding, using Lemma \ref{lem-MomtailF} and $\tau>4$,
    \eqn{
    \expec\big[\indic{a\tilde W_n^* t^2>1}(1+a\tilde W_n^* t^2/2)\big]
    \leq \frac 32 a t^2 \expec\big[\indic{a\tilde W_n^* t^2>1}\tilde W_n^*\big] =o(t^2),
    }
we finally obtain that, uniformly for $t\in [-\eta,\eta]$, there exists a $b>0$ such that
$|\phi_n(t)|\leq 1-b t^2.$ Thus, there exists a constant $a_2>0$ such that
    \eqn{
    \int_{[-\pi,\pi]} |\phi_n(t)|^{l}\frac{dt}{2\pi}
    \leq (1-\vep)^l +\int_{[-\eta,\eta]} (1-b t^2)^{l}\frac{dt}{2\pi}
    \leq \frac{a_2}{l^{1/2}},
    }
which proves \refeq{sums-Xi-bd} for $\delta=2$ and $\tau>4$.

In order to prove \refeq{sums-Xi-bd} for $\tau\in (3,4)$, for which $\delta=\tau-2<2$, we have to
obtain a sharper upper bound on $|\phi_n(t)|$. For this, we identify
$\phi_n(t)=\Re(\phi_n(t))+\i\Im(\phi_n(t))$,
where
    \eqan{
    \Re(\phi_n(t))&=\expec\big[\cos(\tilde W_n^*\sin(t))\mathrm{e}^{-\tilde W_n^*[1-\cos(t)]}\big],
    \quad
    \Im(\phi_n(t))=\expec\big[\sin(\tilde W_n^*\sin(t))\mathrm{e}^{-\tilde W_n^*[1-\cos(t)]}\big],
    }
so that
    \eqn{
    \lbeq{phi-Re-Im}
    |\phi_n(t)|^2=\Re(\phi_n(t))^2+\Im(\phi_n(t))^2.
    }
We start by upper bounding $|\Im(\phi_n(t))|$, by using that
$|\sin(t)|\leq |t|$ for all $t\in \R$,
so that
    \eqn{
    \lbeq{phi-Im-bd}
    |\Im(\phi_n(t))|\leq |t|\expec[\tilde W_n^*]=|t| (1+o(1)).
    }
Further,
    \eqn{
    \Re(\phi_n(t))
    =1-\expec[1-\cos(\tilde W_n^*\sin(t))]
    +\expec\big[\cos(\tilde W_n^* \sin(t)) [\mathrm{e}^{-\tilde W_n^*[1-\cos(t)]}-1]\big].
    }
By the uniform convergence in \refeq{phi-n-conv} and the fact that, for $\eta>0$ small enough,
$\Re(\phi(t))\geq 0$, we only need to derive an upper bound on $\Re(\phi_n(t))$ rather than
on $|\Re(\phi_n(t))|$. For this, we use the fact that $1-\mathrm{e}^{-x}\leq x$ and $1-\cos(t)\leq t^2/2$,
to bound
    \eqan{
    \lbeq{(4.19)}
    \big|\expec\big[\cos(\tilde W_n^* \sin(t)) [\mathrm{e}^{-\tilde W_n^*[1-\cos(t)]}-1]\big]\big|
    &\leq \expec\big[1-\mathrm{e}^{-\tilde W_n^*[1-\cos(t)]}\big]
    \leq [1-\cos(t)]\expec[\tilde W_n^*] \leq \tilde \nu_nt^2 /2.
    }
Further, using \refeq{1-cos-bd} whenever $\tilde W_n^*|t|\leq 1$, so that
also $\tilde W_n^*|\sin(t)|\leq \tilde W_n^*|t|\leq 1$, and $1-\cos(\tilde W_n^*\sin(t))\geq 0$
otherwise, we obtain
    \eqn{
    \lbeq{Re-bd-a}
    \Re(\phi_n(t))
    \leq 1-a\sin(t)^2\expec\big[\big(\tilde W_n^*\big)^2\indic{\tilde W_n^* |t|\leq 1}\big]
    +\tilde \nu_n t^2/2=1-a t^2 \frac{\expec\big[W_n^3\indic{W_n|t|\leq 1}\big]}{\expec[W_n]}
    +\tilde \nu_n t^2/2.
    }
By Lemma \ref{lem-MomtailF}, we have that
    \eqn{
    \lbeq{Re-bd-b}
    \expec\big[W_n^3\indic{W_n|t|\leq 1}\big]\geq
    C_1\big(|t| \vee n^{-1/(\tau-1)}\big)^{\tau-4}.
    }
Combining \refeq{(4.19)} with \refeq{Re-bd-b}, we obtain
that, uniformly in $|t|\leq \eta$ for some small enough $\eta>0$,
    \eqn{
    \lbeq{Re-bd-c}
    \Re(\phi_n(t))
    \leq
    \begin{cases}
    1-2a_{\rm ub}|t|^{\tau-2} &\text{for }|t|\geq n^{-1/(\tau-1)},\\
    1-2a_{\rm ub}t^2 n^{(4-\tau)/(\tau-1)}&\text{for }|t|\leq n^{-1/(\tau-1)},
    \end{cases}
    }
which, combined with \refeq{phi-Re-Im} and \refeq{phi-Im-bd}, shows that, for
$|t|\leq \eta$ and $\eta>0$ sufficiently small,
    \eqn{
    \lbeq{phi-bd-(3,4)}
    |\phi_n(t)|\leq \begin{cases}
    \mathrm{e}^{-a_{\rm ub} t^{2-\tau}}&\text{for }|t|\geq n^{-1/(\tau-1)},\\
    \mathrm{e}^{-a_{\rm ub}t^2 n^{(4-\tau)/(\tau-1)}}
    &\text{for }|t|\leq n^{-1/(\tau-1)}.
    \end{cases}
    }
Thus, there exists a constant $\tilde a_1>0$ such that
    \eqan{
    \int_{[-\pi,\pi]} |\phi_n(t)|^{l}\frac{dt}{2\pi}
    &\leq (1-\vep)^l +\int_{[-\eta,\eta]} \mathrm{e}^{-la_{\rm ub} |t|^{2-\tau}/2}dt
    +\int_{[-n^{-1/(\tau-1)},n^{-1/(\tau-1)}]} \mathrm{e}^{-la_{\rm ub}t^2 n^{(4-\tau)/(\tau-1)}/2}dt\nn\\
    &\leq \frac{\tilde a_1}{l^{1/(\tau-2)}}+\frac{\tilde a_1 n^{(\tau-4)/2(\tau-1)}}{\sqrt{l}}=
    \tilde a_1\big(l^{-1/\delta}+n^{(\tau-4)/2(\tau-1)}l^{-1/2}\big),
    }
which proves \refeq{sums-Xi-bd} for $\tau\in (3,4)$ and with $\delta=\tau-2$.
\qed
\medskip


\noindent
{\it Proof of Proposition \ref{prop-main}(a).} By \refeq{OD-form-applied} and
Lemma \ref{lem-surv-prob},
    \eqan{
    \lbeq{OD-form-again}
    \prob(|\cluster(\Ver)|\geq k)&\leq c (\vep_n^{1/(\delta-1)}\vee n^{-1/(\tau-1)})
    +\tilde a_1\sum_{l=k}^{\infty} \frac{1}{l^{(\delta+1)/\delta}}
    +\tilde a_1 \sum_{l=k}^{\infty} l^{-3/2}\nn\\
    &\leq c \big(\vep_n\vee n^{(\delta-1)/(\tau-1)}\big)^{1/(\delta-1)}+\frac{\tilde a_1 \delta}{k^{1/\delta}}
    +\big(n^{(\tau-4)/2(\tau-1)}\wedge1\big) k^{-1/2}.
    }
The proof is completed by noting that, for $k\geq n^{\delta/(\delta+1)}
=n^{(\tau-2)/(\tau-1)}$,
    \eqn{
    n^{(\tau-4)/2(\tau-1)}k^{-1/2}
    \leq n^{(\tau-4)/2(\tau-1)} n^{(\tau-2)/2(\tau-1)}
    =n^{-1/(\tau-1)}.
    }
Thus, the last term in \refeq{OD-form-again} can be incorporated into the first term,
for the appropriate choice of $a_1$.
This proves the claim in \refeq{clustertail-UB}.
\qed

\paragraph{An upper bound on the expected cluster size: Proof of Proposition \ref{thm-UBclust-size}.}
We now slightly extend the above computation to prove a
bound on the expected cluster size.
We pick $\Lambda>0$ so large that $\tilde\nu_n=(1+\tilde\vep_n)\nu_n
\leq 1-n^{-(\delta-1)/(\delta+1)}$, which is possible
since $\nu_n-1\leq cn^{-(\tau-3)/(\tau-1)}\leq cn^{-(\delta-1)/(\delta+1)}$ by
Corollary \ref{cor-varphi-n}(b) and \refeq{delta-def}.
Then, by \refeq{SD-cluster} and \refeq{T2-SD-T}, as required
    \eqn{
    \expec[|\cluster(\Ver)|]\leq \expec[T^{\sss(2)}]
    \leq \expec[T]=\frac{1}{1-\tilde\nu_n}
    \leq n^{(\delta-1)/(\delta+1)}.
    }

\section{A lower bound on the cluster tail}
\label{sec-LBcluster_tail}
In this section, we prove a lower bound on the cluster tail.
The key ingredient in the proof of Proposition \ref{prop-main}(b) is
again the coupling to branching processes. Note the explicit coupling
between the cluster size $|\cluster(\Ver)|$ and the total progeny $T^{\sss(2)}$
described there. We can then bound
    \eqn{
    \lbeq{split-cluster-tail-lb}
    \prob(|\cluster(\Ver)|\geq k)
    \geq \prob(T^{\sss(2)}\geq 2k, |\cluster(\Ver)|\geq k)
    =\prob(T^{\sss(2)}\geq 2k) -\prob(T^{\sss(2)}\geq 2k, |\cluster(\Ver)|<k).
    }
The following lemmas contain bounds on both contributions:

\begin{lemma}[Lower bound tail total progeny]
\label{lem-T-lbd}
Under the assumptions of Proposition \ref{prop-main}(b),
there exists a constant $a_2>0$ such that,
for all $k\leq \vep n^{(\delta-1)/(\delta+1)}$,
    \eqn{
    \lbeq{T-lbd}
    \prob(T^{\sss(2)}\geq k)\geq \frac{2a_2}{k^{1/\delta}}.
    }
\end{lemma}

\begin{lemma}[Upper bound cluster tail coupling]
\label{lem-coupling-bd}
Fix $\vep>0$ sufficiently small.
Under the assumptions of Proposition \ref{prop-main}(b), there exists constants $c,q>0$ such that,
for all $k\leq \vep n^{(\delta-1)/(\delta+1)}$,
    \eqn{
    \lbeq{coupling-bd}
    \prob(T^{\sss(2)}\geq 2k, |\cluster(\Ver)|<k)\leq \frac{c\vep^{q}}{k^{1/\delta}}.
    }
\end{lemma}
\medskip

\noindent
{\it Proof of Proposition \ref{prop-main}(b) subject to
Lemmas \ref{lem-T-lbd}-\ref{lem-coupling-bd}.} Recall \refeq{split-cluster-tail-lb}, and substitute the
bounds in Lemmas \ref{lem-T-lbd}-\ref{lem-coupling-bd} to conclude that
    \eqn{
    \prob(|\cluster(\Ver)|\geq k)\geq \frac{2a_2}{(2k)^{1/\delta}}-\frac{c\vep^{p}}{k^{1/\delta}}
    \geq \frac{a_2}{k^{1/\delta}},
    }
when $\vep>0$ is so small that $2^{1-1/\delta} a_2-c\vep^{p}\geq a_2$. This is possible, since $\delta>1$.
\qed

\medskip

\noindent
{\it Proof of Lemma \ref{lem-T-lbd}.} We start by noting that
    \eqn{
    \lbeq{lb-TP-BP}
    \prob(T^{\sss(2)}\geq k)\geq \prob(T^{\sss(2)}\geq k, X_1=1)=\prob(T\geq k-1)\prob(X_1=1)
    \geq \prob(T\geq k)\prob(X_1=1).
    }
Note that $\prob(X_1=1)=\expec[W_n\e^{-W_n}]=\expec[W\mathrm{e}^{-W}]+o(1),$
which remains strictly positive. Thus, it suffices to prove a lower bound on $\prob(T\geq k)$.
For this, we bound
    \eqan{
    \prob(T\geq k)
    &\geq\sum_{l=k}^{\infty} \prob(T=l)=\sum_{l=k}^{\infty} \frac{1}{l} \prob \big(\sum_{i=1}^l X_i = l-1\big)
    \geq \sum_{l=k}^{2k} \frac{1}{l} \prob \big(\sum_{i=1}^l X_i = l-1\big).
    }

We prove the bounds for $\tau\in (3,4)$ and $\tau>4$ simultaneously,
the latter being somewhat simpler. We shall follow a large
part of the analysis for the upper bound in Proposition \ref{prop-sumsXi}.
Recall \refeq{FIT-applied}, to get
    \eqn{
    \lbeq{FIT-applied-rep}
    \prob \big(\sum_{i=1}^l X_i = l-1\big)=\int_{[-\pi,\pi]} \Re\Big(\mathrm{e}^{\i t} \big(\mathrm{e}^{-\i t}\phi_n(t)\big)^{l}\Big)\frac{dt}{2\pi}.
    }
It is not hard to see that, by the arguments in the proof of Proposition \ref{prop-sumsXi}
(in particular, recall \refeq{phi-bd-(3,4)}),
    \eqn{
    \int_{[-\pi,\pi]} \Big|(\mathrm{e}^{\i t}-1) \big(\mathrm{e}^{-\i t}\phi_n(t)\big)^{l}\Big|\frac{dt}{2\pi}
    \leq cl^{-2/\delta},
    }
and, with $K_l=Kl^{-1/\delta}$,
    \eqn{
    \int_{[-\pi,\pi]\setminus [-K_l,K_l]} \big(\mathrm{e}^{-\i t}\phi_n(t)\big)^{l}\frac{dt}{2\pi}
    \geq -{\mathrm e}^{-cK^{\tau-2}}l^{-1/\delta},
    }
so we are left to prove a lower bound of the form $cl^{-1/\delta}$ for $\int_{[-K_l,K_l]} \Re(\big(\mathrm{e}^{-\i t}\phi_n(t)\big)^{l})\frac{dt}{2\pi}.$
Note that, by the uniform convergence of $\mathrm{e}^{-\i t}\phi_n(t)$ to
$\mathrm{e}^{-\i t}\phi(t)$, we have that $\mathrm{e}^{-\i t}\phi_n(t)\rightarrow 1$
uniformly for $t\in [-K_l,K_l]$.

Let $\varphi_n(t)=\mathrm{e}^{-\i t}\phi_n(t), \varphi(t)=\mathrm{e}^{-\i t}\phi(t)$.
By scaling,
    \eqn{
    \lbeq{grand-aim}
    \int_{[-K_l,K_l]} \Re\big(\varphi_n(t)^{l}\big)\frac{dt}{2\pi}
    =l^{-1/\delta}\int_{[-K,K]} \Re\big(\varphi_n(tl^{-1/\delta})^{l}\big)\frac{dt}{2\pi}.
    }
We rewrite
    \eqn{
    \varphi_n(t)^{l}
    =\varphi(t)^{l}\Big(1+\frac{\varphi_n(t)-\varphi(t)}{\varphi(t)}\Big)^{l},
    }
so that
    \eqn{
    \Re\big(\varphi_n(t)^{l}\big)
    =\Re\Big(\varphi(t)^{l}\Big) \Re\Big(1+\frac{\varphi_n(t)-\varphi(t)}{\varphi(t)}\Big)^{l}
    -\Im\Big(\varphi(t)^{l}\Big) \Im\Big(1+\frac{\varphi_n(t)-\varphi(t)}{\varphi(t)}\Big)^{l}.
    }
When $t\rightarrow 0$, with $\alpha=(\tau\wedge 4)-2$,
    \eqn{
    \varphi(t)={\mathrm e}^{-c|t|^{\alpha}(1+o(1))},
    }
since $W^*-1$ is in the domain of attraction of an $\alpha$-stable distribution when
\refeq{F-bound-tau>4} or \refeq{F-bound-tau(3,4)b} hold, and thus also $X=\Poi(W^*)-1$ is.
Thus, pointwise in $t$ as $l\rightarrow \infty$,
    \eqn{
    \Re\Big(\varphi(tl^{-1/(\tau-2)})^{l}\Big)=
    {\mathrm e}^{-c|t|^{\alpha}}(1+o(1)),
    \qquad
    \Im\Big(\varphi(tl^{-1/\delta})^{l}\Big)=o(1).
    }
Therefore, also using the fact that
    \eqn{
    \Big|\Im\Big(\varphi(tl^{-1/\delta})^{l}\Big) \Im\Big(1+\frac{\varphi_n(tl^{-1/\delta})-\varphi(tl^{-1/\delta})}{\varphi(tl^{-1/\delta})}\Big)^{l}\Big|
    \leq  \big|\varphi_n(tl^{-1/\delta})^l\big| \leq {\mathrm e}^{-c |t|^{\alpha}},
    }
which is integrable, dominated convergence gives that,
for every $K>0$,
    \eqn{
    \int_{[-K,K]} \Im\Big(\varphi(tl^{-1/\delta})^{l}\Big) \Im\Big(1+\frac{\varphi_n(tl^{-1/\delta})-\varphi(tl^{-1/\delta})}{\varphi(tl^{-1/\delta})}\Big)^{l}\frac{dt}{2\pi}
    =o(1).
    }
By Corollary \ref{cor-varphi-n}(a) and the fact that $|\varphi(t)|\geq 1/2$
for all $|t|\leq K_l$,
    \eqan{
    \Re\Big(1+\frac{\varphi_n(t)-\varphi(t)}{\varphi(t)}\Big)^{l}
    &\geq (1-2|\varphi_n(t)-\varphi(t)|)^l\geq 1-2l|\varphi_n(t)-\varphi(t)|\nn\\
    &\geq 1-2l cn^{-(\tau-2)/(\tau-1)}+2cl|t|\big(|\vep_n|+n^{-(\tau-3)/(\tau-1)}\big)\nn\\
    &\geq 1-2c (ln^{-(\tau-2)/(\tau-1)}) -2c|t|l^{1/(\tau-2)} (ln^{-(\tau-2)/(\tau-1)})^{(\tau-3)/(\tau-2)}.
    \lbeq{Re-bd-a}
    }
Now we use the fact that $l\leq 2k\leq 2\vep n^{\delta/(\delta+1)}$,
where $\delta/(\delta+1)=(\tau-2)/(\tau-1)$ for $\tau\in (3,4)$ and
$\delta/(\delta+1)=2/3<(\tau-2)/(\tau-1)$ for $\tau>4$,
while $|\vep_n|\leq \Lambda n^{-(\delta-1)/(\delta+1)}$,
where $(\delta-1)/(\delta+1)=(\tau-3)/(\tau-1)$ for $\tau\in (3,4)$ and
$(\delta-1)/(\delta+1)=1/3<(\tau-3)/(\tau-1)$ for $\tau>4$. Therefore,
for $\tau>4$, the left-hand side of \refeq{Re-bd-a} is $1+o(1)$,
while, for $\tau\in (3,4)$ we can bound it as
    \eqan{
    \lbeq{Re-bd}
    \Re\Big(1+\frac{\varphi_n(t)-\varphi(t)}{\varphi(t)}\Big)^{l}
    &\geq 1-4c\vep -4c\vep^{(\tau-3)/(\tau-2)} (1+\Lambda)|t|l^{1/(\tau-2)}\nn\\
    &\geq 1-4c\vep -4c\vep^{(\tau-3)/(\tau-2)} (1+\Lambda)K\geq 1/2,
    }
when $4c\vep(1+\Lambda)+4c\vep^{(\tau-3)/(\tau-2)}K\leq 1/2$. In particular,
for all $|t|\leq K_l$, the left-hand side of \refeq{Re-bd} is non-negative.
Therefore, we arrive at the claim that, for $\vep>0$ sufficiently small,
there exists $C=C(\vep, \Lambda)>0$ such that, as $l\rightarrow \infty$,
    \eqn{
    \int_{[-K,K]} \Re\big(\varphi_n(tl^{-1/\delta})^{l}\big)\frac{dt}{2\pi}
    \geq C\int_{[-K,K]} {\mathrm e}^{-c|t|^{\alpha}}dt
    +o(1).
    }
This completes the proof of Lemma \ref{lem-T-lbd}.
\qed
\medskip

\noindent
{\it Proof of Lemma \ref{lem-coupling-bd}.}
In the proof, we will make repeated use of the following lemma:

\begin{lemma}[Upper bound on tail probabilities of random sums]
\label{lem-random-sums-bd}
Suppose that $(T_l)_{l\geq 1}$ are i.i.d.\ random variables
for which there exist constants $K>0$ and $\delta>1$ such that
    \eqn{
    \lbeq{ass-random-sums-bd}
    \prob(T_l\geq k)\leq  K k^{-1/\delta},
    }
and let $M$ be independent from $(T_l)_{l\geq 1}$. Then,
there exists a constant $c_{\delta}>0$ such that
    \eqn{
    \prob\Big(\sum_{l=1}^M T_l\geq k\Big)\leq  c_{\delta} K \expec[M] k^{-1/\delta}.
    }
\end{lemma}

\proof We split, depending on whether there exists a $j\in [M]$ such that $T_j\geq k$ or not,
    \eqan{
    \prob(\sum_{l=1}^M T_l\geq k)&\leq \prob(\exists j\in [M]: T_j\geq k)
    +\prob\Big(\sum_{l=1}^M T_l\indic{T_l\leq k}\geq k\Big)\nn\\
    &\leq \expec[M] \prob(T_1\geq k) + \frac 1k \expec\Big[\sum_{l=1}^M T_l \indic{T_l\leq k}\Big]\nn\\
    &\leq K \expec[M] k^{-1/\delta}+ \frac 1k \expec[M]\sum_{l=1}^k \prob(T_1\geq l)\nn\\
    &\leq K \expec[M] \Big(k^{-1/\delta}+\frac1k\sum_{l=1}^k l^{-1/\delta}\Big)
    \leq c_{\delta} K \expec[M] k^{-1/\delta},
    }
where in the second inequality,
we use Boole's inequality for the first term and Markov's inequality for the second.
\qed

\noindent
We write
    \eqn{
    \prob(T^{\sss(2)}\geq 2k, |\cluster(\Ver)|<k)
    =\sum_{t=1}^{k-1} \prob(T^{\sss(2)}\geq 2k, |\cluster(\Ver)|=t).
    }
For $j\in [n]$, let $M_j(t)$ denote the number of times the mark
$j$ is drawn in the first $t$ draws, and define $N_t=\{j\in [n]\colon M_j(t)\geq 2\}$.
For $j\in [n]$, let $(T_{j,s})_{s=1}^{\infty}$ be an i.i.d.\ sequence
of random variables where $T_{j,1}$ is the total progeny of a branching process which has
a $\Poi(w_j)$ offspring in the first generation, and offspring
distribution $\Poi(W_n^*)$ in all later generations.
When $|\cluster(\Ver)|=t$, but $T^{\sss(2)}\geq 2k>t$, then we must have that
the total progeny of the thinned vertices is at least $k$. A vertex is thinned
precisely when its mark is chosen at least twice, so that we can bound
    \eqn{
    \prob\big(T^{\sss(2)}\geq 2k, |\cluster(\Ver)|=t\big)
    \leq \prob\Big(T^{\sss(2)}\geq 2k, |\cluster(\Ver)|=t, \sum_{j\in N_t} \tilde T_{j,t}\geq k\Big),
    }
where, since all repetitions of vertex $j$ after the first are thinned,
    \eqn{
    \lbeq{Tjt-def}
    \tilde T_{j,t}=\sum_{s=1}^{M_j(t)-1} T_{j,s}.
    }
Since $t\mapsto N_t$ and $t\mapsto M_j(t)$ are non-decreasing, we arrive at
    \eqn{
    \prob\big(T^{\sss(2)}\geq 2k, |\cluster(\Ver)|<k\big)
    \leq \prob\Big(T^{\sss(2)}\geq 2k,\sum_{j\in N_k} \tilde T_{j,k}\geq k\Big).
    }
Now, as in the proof of Lemma \ref{lem-random-sums-bd},
we split the event $\sum_{j\in N_k} \tilde T_{j,k}\geq k$ into
the event where all $\tilde T_{j,k}\leq k$ and the event where there exists a
$j\in N_k$ such that $\tilde T_{j,k}\geq k$, to arrive at
    \eqan{
    \lbeq{clever-split}
    \prob\big(T^{\sss(2)}\geq 2k, |\cluster(\Ver)|<k\big)
    &\leq \sum_{j=1}^n \prob(T^{\sss(2)}\geq 2k,j\in N_k,\tilde T_{j,k}\geq k)
    +\prob(T^{\sss(2)}\geq 2k,\sum_{j\in N_k} \tilde T_{j,k}\indic{\tilde T_{j,k}\leq k}\geq k)\nn\\
    &\leq\sum_{j=1}^n \Big[\prob\big(T^{\sss(2)}\geq 2k,j\in N_k,\tilde T_{j,k}\geq k\big)
    +\frac 1k \sum_{l=1}^k \prob\big(T^{\sss(2)}\geq 2k,j\in N_k,\tilde T_{j,k}\geq l\big)\Big],\nn
    }
the last bound by the Markov inequality.

Let $(K_1, K_2)$ be the first two times before $\tildeT$ for which $M_{K_1}=M_{K_2}=j$. Then,
noting that $T^{\sss(2)}\geq k_1$ needs to occur when $K_1=k_1$,
    \eqan{
    \prob\big(T^{\sss(2)}\geq 2k,j\in N_k,\tilde T_{j,k}\geq l\big)
    &\leq \sum_{1\leq k_1<k_2\leq k}\prob(T^{\sss(2)}\geq k_1, (K_1,K_2)=(k_1,k_2))\nn\\
    &\qquad \times \prob\Big(\tilde T_{j,k}\geq l\mid (K_1,K_2)=(k_1,k_2)\Big).
    }
We wish to apply Lemma \ref{lem-random-sums-bd} to $\tilde T_{j,k}$ in \refeq{Tjt-def},
where $M$ is $M_j(k)-1$ conditioned on $(K_1,K_2)=(k_1,k_2)$. When $k\leq \vep n^{(\tau-2)/(\tau-1)}$,
    \eqn{
    \expec[M_j(k)-1\mid (K_1,K_2)=(k_1,k_2)]
    =2+(k-k_2)\frac{w_j}{\ell_n}\leq 2+\vep n^{(\tau-2)/(\tau-1)}\frac{w_1}{\ell_n}=O(1).
    }
In order to apply Lemma \ref{lem-random-sums-bd}, we proceed
by checking \refeq{ass-random-sums-bd}.
Note that $(T_{j,s})_{s\geq 1}$ is i.i.d.\ with
$T_{j,1}=1+\sum_{l=1}^{P_j} T_l$, where $P_j\sim \Poi(w_j)$ and where
$(T_l)_{l\geq 1}$ is an i.i.d.\ sequence of total progenies of
branching processes with offspring distribution $\Poi(W_n^*)$.
Thus, by Lemma \ref{lem-random-sums-bd}, for which the assumption follows from
Proposition \ref{prop-main}(a) (using the fact that $\vep_n\leq 0$), we have
    \eqn{
    \lbeq{Tj1-tail-bd}
    \prob(T_{j,1}\geq l)
    =\prob\Big(\sum_{l=1}^{P_j} T_j \geq l-1\Big)
    \leq c_\delta a_1 \expec[P_j]l^{-1/\delta}= \frac{c_\delta a_1 w_j}{l^{1/\delta}}.
    }
Therefore, applying Lemma \ref{lem-random-sums-bd} yields that
    \eqan{
    \lbeq{Tj1-tail-bd-2}
    \prob\big(T^{\sss(2)}\geq 2k,j\in N_k,\tilde T_{j,k}\geq l\big)
    &\leq \sum_{1\leq k_1<k_2\leq k}\prob(T^{\sss(2)}\geq k_1, (K_1,K_2)=(k_1,k_2))
    \frac{c_\delta a_1 w_j}{l^{1/\delta}}.
    }
The factor $w_j$ appearing in \refeq{Tj1-tail-bd-2} can be harmful
when $w_j$ is large. To resolve this problem, for $j$ such
that $w_j\geq k^{1/\delta}$, we simply bound the restrictions
on $\tilde T_{j,k}$ away, so that
    \eqan{
    \lbeq{clever-split-rep}
    &\prob\big(T^{\sss(2)}\geq 2k, |\cluster(\Ver)|<k\big)\\
    &\leq \quad 2\sum_{j=1}^n \sum_{1\leq k_1<k_2\leq k}
    \indic{w_j<k^{1/\delta}}
    \prob\big(T^{\sss(2)}\geq k_1, (K_1,K_2)=(k_1,k_2)\big)\nn\\
    &\qquad \qquad+\sum_{j=1}^n \sum_{1\leq k_1<k_2\leq k}
    \indic{w_j<k^{1/\delta}}c_\delta a_1 w_j
    \prob\big(T^{\sss(2)}\geq k_1, (K_1,K_2)=(k_1,k_2)\big)
    \Big(\frac{1}{l^{1/\delta}}+\frac 1k \sum_{l=1}^k\frac{1}{l^{1/\delta}}\Big).\nn
    }
Performing the sum over $l$ and bounding $\prob((K_1,K_2)=(k_1,k_2))\leq (w_j/\ell_n)^2$
leads to
    \eqan{
    \prob\big(T^{\sss(2)}\geq 2k, |\cluster(\Ver)|<k\big)
    &\leq cc_\delta a_1 \sum_{j=1}^n \sum_{1\leq k_1<k_2\leq k} \prob\big(T^{\sss(2)}\geq k_1\big)
    \big(\frac{w_j}{k^{1/\delta}}\wedge 1\big)\big(\frac{w_j}{\ell_n}\big)^2\\
    &\leq cc_\delta a_1 \sum_{j=1}^n \sum_{1\leq k_1<k_2\leq k} \big(\frac{w_j}{\ell_n}\big)^2 a_2 k_1^{-1/\delta}\big(\frac{w_j}{k^{1/\delta}}\wedge 1\big)\nn\\
    &\leq cc_\delta a_1 \frac{nk^{2-1/\delta}}{\ell_n^2} \expec\Big[W_n\big(\frac{W_n^2}{k^{1/\delta}}\wedge 1\big)\Big].\nn
    }

When $\tau>4$, we have that $\delta=2$ and $k\leq \vep n^{2/3}$, so that we can bound the above by
    \eqn{
    \prob\big(T^{\sss(2)}\geq 2k, |\cluster(\Ver)|<k\big)
    \leq C_{\delta} \frac{k^{2-2/\delta}}{n}\expec[(W^*)^2]
    = C_{\delta} \frac{k}{n}=C_{\delta} \frac{k^{3/2}}{n} k^{-1/2}\leq C_{\delta} \vep^{3/2} k^{-1/2},
    }
for some constant $C_\delta>0$, so that \refeq{coupling-bd} follows with $q=3/2$.

When $\tau\in (3,4)$, we use Lemma \ref{lem-MomtailF},
now with $k\leq \vep n^{(\tau-2)/(\tau-1)}$ and
$\delta=\tau-2$, to obtain
    \eqan{
    \prob\big(T^{\sss(2)}\geq 2k, |\cluster(\Ver)|<k\big)
    &\leq C_\delta\frac{k^{2-2/\delta}}{n} \expec\Big[W_n^3
    \indic{W_n\leq k^{1/\delta}}\Big]
    +C_\delta\frac{k^{2-1/\delta}}{n} \expec\Big[W_n^2\indic{W_n>k^{1/\delta}}\Big]\nn\\
    &\leq C_\delta \Big(\frac{k^{(\tau-1)/(\tau-2)}}{n}\Big) k^{-1/(\tau-2)}\leq C_\delta \vep^{(\tau-1)/(\tau-2)} k^{-1/(\tau-2)},
    }
for some constant $C_\delta>0$, so that \refeq{coupling-bd} follows with $q=(\tau-1)/(\tau-2)>1.$
\qed

%

\paragraph{Acknowledgements.}
This work was supported in part by the Netherlands
Organisation for Scientific Research (NWO). The author gratefully
acknowledges the hospitality of the Mittag-Leffler Institute (MLI),
where part of this work was carried out, and the lively discussions
there with Svante Janson, Tomasz {\L}uczak, Michal Karo\'nski and
Tatyana Turova, as well as Johan van Leeuwaarden and two anonymous referees
for useful comments on a preliminary version that greatly improved the
presentation.

\bibliographystyle{plain}
\bibliography{../../onderzoek/bib/bib}

\end{document}